\documentclass[11pt]{article}
\newtheorem{theorem}{Theorem}[section]
\newtheorem{lemma}{Lemma}[section]

\newtheorem{corollary}{Corollary}[section]

\newcommand{\R}{\mbox{\rm I$\!$R}}

\newcommand{\Z}{\mbox{\rm \lower0.3pt\hbox{$\angle\!\!\!$}Z}}
\newcommand{\C}{\mbox{\rm $~\vrule height6.5pt width0.5pt depth0.3pt\!\!$C}}
\newcommand{\Hs}{{\cal H}}
\newcommand{\li}{\infty}

\newcommand{\lb}{\left (}
\newcommand{\rb}{\right )}
\newcommand{\la}{\langle }
\newcommand{\ra}{\rangle }

\newcommand{\ddem}{\ind {\it Proof: }}
\newcommand{\ok}{\bullet}
\newcommand{\edem}{. $\ok$}

\newcommand{\ind}{\hspace*{0.3cm} }

\newcommand{\dis}{\displaystyle}
\newcommand{\bc}{\begin{center}}
\newcommand{\ec}{\end{center}}
\newcommand{\benq}{\begin{eqnarray}}
\newcommand{\be}{\begin{equation}}

\newcommand{\ee}{\end{equation}}
\newcommand{\eenq}{\end{eqnarray}}
\newcommand{\ba}{\begin{array}{l}}
\newcommand{\ea}{\end{array}}

\newcommand{\bfg}{\begin{figure}[h]}
\newcommand{\efg}{\end{figure}}

\newcommand{\lp}{\left [ }
\newcommand{\rp}{\right ] }
\newcommand{\ov}{\overline}

\newcommand{\es}{\hskip 0.3cm}

\newcommand{\al}{\alpha}
\newcommand{\ld}{\lambda}

\newcommand{\eps}{\epsilon}

\newtheorem{proposition}{Proposition}[section]

\newcommand{\argmin}{{\rm argmin}}

\newcommand{\der}{\partial}
\newcommand{\uds}{\ddot{u}}
\newcommand{\udp}{\dot{u}}
\newcommand{\gam}{\gamma}
\usepackage{subeqn}
\begin{document}
\begin{center}
 {\bf \Large Asymptotic behavior of  second-order dissipative evolution equations combining
  potential with non-potential effects }
\end{center}

\begin{center}
 {\bf Hedy Attouch}\footnote{ Institut de Math\'ematiques et de Mod\'elisation de
        Montpellier, UMR CNRS 5149, CC 51, Universit\'e Montpellier II,
        Place Eug\`ene Bataillon, 34095 Montpellier cedex 5, France.
            ({\tt  attouch@math.univ-montp2.fr}). Supported by French ANR
grant ANR-08-BLAN-0294-03.} \ and \
       {\bf Paul-Emile Maing\'e}\footnote{Universit\'e des Antilles-Guyane,
       D.S.I., CEREGMIA,  Campus de Schoelcher, 97233 Cedex,
       Martinique, F.W.I.
            ({\tt  Paul-Emile.Mainge@martinique.univ-ag.fr}).} \\
\end{center}

\begin{center}
\textit{Dedicated to Alain Haraux on the occasion of his 60th birthday}
\end{center}

\bigskip
\bigskip

\noindent{\bf Abstract} {\small   We investigate,  in the setting of
a real Hilbert space ${\cal H}$, the asymptotic behavior, as time
$t$ goes to infinity, of trajectories of second-order evolution
equations
 $$ \uds (t) + \gam \udp(t)  +  \nabla \phi(u(t)) + A(u(t)) = 0,$$

\noindent where  $\nabla \phi$ is the gradient operator of a convex
differentiable potential function $\phi: \Hs \to \R$, $A: \Hs \to
\Hs$ is a   maximal monotone operator which is assumed to be
$\lambda$-cocoercive, and $\gamma > 0$ is a  damping parameter.
Potential and non-potential effects are associated  respectively to
$\nabla \phi$ and $A$. We prove that, under condition   $\
\lambda{\gam}^2 > 1,  \ $  each trajectory asymptotically  weakly
converges to a zero of $ \nabla \phi + A $. This condition, which
only  involves the non-potential  operator  and the damping
parameter, is  sharp and consistent with time rescaling. Passing
from  weak to  strong convergence of the trajectories is obtained by
introducing  an asymptotically vanishing Tikhonov-like regularizing
term. As special cases,  we recover  the asymptotic analysis of  the
heavy ball with friction dynamic attached to a convex potential, the
second-order gradient-projection dynamic, and the second-order
dynamic governed by the Yosida  approximation of a general maximal
monotone operator. The breadth and flexibility
 of the proposed framework is  illustrated through  applications in the areas of  constrained optimization,
dynamical approach to Nash equilibria for noncooperative games, and asymptotic stabilization in the case of a continuum of equilbria. }

\bigskip

\bigskip


\bigskip


 \pagestyle{myheadings} \thispagestyle{plain} \markboth{H. ATTOUCH, P.E.
MAINGE}{Asymptotic behavior of  second-order dissipative evolution
equations}
\section{\large Introduction}
Throughout this paper, $\Hs$ is a real Hilbert space.  Its scalar
product is denoted by $\la ,\ra$ and the associated norm by  $|
\cdot  |$. We wish to investigate  the asymptotic behavior, as the
time variable $t$ goes to $+\infty$, of  second order evolution
equations

 \be \label{pb000}
\uds (t) + \gam \udp(t)  +   \mathcal A(u(t)) = 0
\ee

\noindent   where $\mathcal A: \Hs \to \Hs $ is a maximal monotone
operator and $\gam >0$ is a positive damping parameter. Here and henceforth, a dot $ ( \ \dot{ } \ ) $  denotes first-order differentiation with respect to time $t$, and  a double dot
 $ ( \ \ddot{ } \ )$ denotes
second-order differentiation.
 In order to grasp the respective influence of potential and non-potential operators in such inertial systems,
 and their effects on the asymptotic behavior of  trajectories,
  we  consider the class of maximal monotone operators which can be splitted up into the sum of two terms  $\mathcal A   =  \nabla \phi + A$
 where
\begin{itemize}
\item  $ \nabla \phi : \Hs \to \Hs $  is the gradient operator of  a  convex and continuously differentiable function $\phi :\Hs \to \R$;
\item  $A : \Hs \to \Hs$ is   a maximal monotone operator that is
assumed to be cocoercive, which means that there exists some  constant  $\lambda > 0$ such that
\be \label{cc}
\forall (v,w) \in \Hs^2  \mbox{  }   \mbox{  }  \la Av - Aw , v - w \ra  \geq \lambda |Av-Aw|^2 .\es
\ee
\end{itemize}
An operator $A : \Hs \to \Hs$ which satisfies (\ref{cc})  is said to
be  $\lambda$-cocoercive (the relevancy of this hypothesis with
respect to applications is examined below). Note that $A$
$\ld$-cocoercive implies that $A$ is $(1/\ld)$-Lipschitz continuous.
\smallskip

An outline of the present  work  is as follows:

1. In section 2, we study the asymptotic behavior  of second-order
autonomous evolution systems governed by such operators $\mathcal A
=  \nabla \phi + A$, namely

\be \label{pb00} \uds (t) +\gam \udp(t)   + \nabla \phi(u(t))  +
A(u(t)) = 0, \ee

\noindent with $\phi$ convex and continuously differentiable, $A$ \ $\lambda$-cocoercive for some $\lambda>0$,
and  $\gam > 0$ as  a  damping parameter.
Let us denote by

\be \label{equil}
 S:=\{ v \in \Hs \es | \es  \nabla \phi(v) + Av = 0\}
\ee

\noindent the set of equilibria and suppose that $S\neq \emptyset$.
In Theorem \ref{wac} we establish   that, under the sole assumption

\be \label{cf}
\lambda{\gam}^2 > 1,
\ee

\noindent each trajectory $t \to u(t)$ of (\ref{pb00})
 weakly converges in $\Hs$,  as $t \to
+\infty$,  to an element of $S$.
This condition, which only  involves the non-potential part of the maximal monotone operator governing the equation,
is proved to be sharp, in the sense  that  one can exhibit situations where condition  $\ld \gam^2 < 1$ does
not ensure convergence of all trajectories.

When the operator $\nabla \phi + A$ is
 strongly monotone we prove  strong asymptotic convergence of the trajectories.

In the  potential case, $\mathcal A   =  \nabla \phi $
(corresponding to  $A = 0$), taking advantage of the fact that, in
the above result, no restrictive assumption is made on the potential operator
$\nabla \phi$, we recover the asymptotic convergence result for the
so-called heavy ball with friction dynamical system (Alvarez
\cite{alv})

\be \label{pb01}
\uds (t) +\gam \udp(t)  +  \nabla \phi(u(t))  = 0.
\ee

\noindent In recent years, because of its rich connections with mechanics and optimization, this system has been the object of active research,
 see \cite{alv},\cite{ant},  \cite{acb}, \cite{agr},  and the references therein.

In the non-potential case, $\nabla \phi = 0$, (\ref{pb00}) becomes

\be \label{pb02}
\uds (t) +\gam \udp(t)  +  A(u(t))  = 0
\ee

\noindent with $A$ a cocoercive operator. Equation (\ref{pb02}) covers  several situations of practical interest:

\begin{itemize}
   \item  $A = I - T$ where $T : \Hs \to \Hs$ is a
contraction. It can be  easily checked  that $A$ is
$(1/2)$-cocoercive. Condition (\ref{cf}) gives $\gam > \sqrt{2}$.
When specialized to this situation, Theorem \ref{wac} yields
convergence results for the second order gradient-projection
dynamical system, first established in \cite{alat2} and \cite{ant}.
   \item
   $A = B_{\lambda}$ where  $B_{\lambda}$ (with parameter $\ld >0$)
   is the Yosida approximation  of  a general maximal monotone operator
   $B: \Hs \to 2^{\Hs}$, (see \cite{Bre}). One can easily verify that
    $ B_{ \lambda}$ is
   $\ld$-cocoercive. Noticing that $B$ and
$B_{\lambda}$ have the same zeroes, we shall derive new inertial second-order dynamical approach to the set of zeroes $B^{-1} (0)$, for $B$ a general maximal monotone operator.
\end{itemize}

\smallskip

2. In section 3,  we introduce a  Tikhonov-like regularizing
term $\nabla \Theta (u(t))$ with  vanishing coefficient  $\eps(t)$
in the above dynamic, and consider the nonautonomous system

\be \label{tik} \uds (t) +\gam \udp(t)   + \nabla \phi(u(t))+
A(u(t)) +  \eps(t) \nabla \Theta (u(t)) = 0, \ee

\noindent where the function $\Theta: {\cal H} \to \R $ is supposed
to be convex, differentiable, with  $\nabla \Theta$ Lipschitz
continuous
  and strongly monotone, while
  $\eps :[0, +\li)  \to [0,+\li)$ is a function of class $C^1$ such that  $\eps (t) \to 0^+$ as $t
\to +\li$.

In Theorem \ref{imp}, we   establish that, under condition
(\ref{cf}),  and the slow vanishing condition on $\epsilon(\cdot)$

\be \label{sleps} \dis \int_0^{+\li} \eps(s) ds = +\li, \ee

\noindent  each trajectory of (\ref{tik}) strongly converges as $t
\to +\li$ to $u_*$, which is  the unique minimizer of $\Theta$ over
the set $S := (\nabla \phi + A)^{-1}(0)$. Equivalently, $u_*$ solves
the following variational inequality problem: find $u_* \in S$ such
that
 \be \label{pbv} \la \nabla \Theta (u_*), v-u_* \ra
 \ge 0 \es \forall v\in S.
 \ee

Note that the above regularization technique allows both to obtain
strong convergence of the trajectories,
 together with a limit which no longer  depends  on the initial data.
 As a consequence of Theorem \ref{imp}, we recover various results which have been devoted to
  the Tikhonov dynamics, see \cite{acz}, \cite{copeso} and the references
   therein.

\smallskip

3. In the final section, we briefly outline  situations where  our results can be applied.
 Indeed, the convergence results obtained in Theorems  \ref{wac} and \ref{imp} offer promising views on numerical  optimization and algorithms (by time discretization as in   \cite{alv}, \cite{alat}, \cite{atcom}, \cite{combh}), on the modeling of inertial dynamical approach to Nash equilibria in decision sciences and game theory (see \cite{as1}, \cite{copeso}), and on dissipative dynamical systems and PDE.'s (as in \cite{alat3} for the damped wave equation).

%
%
 \section{\large Asymptotic convergence results}
This section is devoted to the study of the asymptotic behavior, as  time variable $t$ goes to infinity,
 of   trajectories of  (\ref{pb00}).

\subsection{Weak asymptotic convergence results}

\noindent Throughout this section we make the following assumptions:

\smallskip

(H1)  \begin{minipage}[t]{12.cm} $\phi : \Hs \to \R$ is a
 convex differentiable function whose gradient $\nabla \phi$ is
 Lipschitz continuous on the bounded subsets of $\Hs$;
 \end{minipage}

\smallskip

(H2) \begin{minipage}[t]{12.cm} $A: \Hs \to \Hs$ is maximal monotone
and $\ld$-cocoercive for some $\ld >0$.
\end{minipage}

\medskip

\noindent Note that  the cocoerciveness of  $A$ implies that $A$ is
Lipschitz continuous, so that $ \nabla \phi + A$ is Lipschitz
continuous on bounded sets. By Cauchy-Lipschitz theorem,  for any
initial  data $u(0)=u_0$,  $\dot{u}(0)=v_0$ with $(u_0,v_0)$ in
$\Hs^2$ there exists a unique local solution to the Cauchy problem
\be \label{pbn} \label{pb0} \left [ \ba
  \uds (t) +\gam \udp(t)   + \nabla \phi(u(t)) +   A(u(t)) = 0,\\
  \\
 u(0)=u_0, \es \dot{u}(0)=v_0. \ea \right.
 \ee

\noindent The fact that $u$ is infinitely extendible to the right
follows from a uniform bound on $|\udp(t)|$ as given in the proof of
Theorem \ref{wac}, which provides a unique classical global solution
$u \in C^2([0,+\li); \Hs)$ of (\ref{pbn}).
\medskip

\noindent The following theorem establishes the weak asymptotic
convergence property of $u(\cdot),$ solution of  (\ref{pb0}),
  under the sole assumption
 $\lambda{\gam}^2 > 1$.

\begin{theorem} \label{wac}
Let us suppose that   {\rm (H1)-(H2) } hold  with  $S := ( \nabla
\phi + A)^{-1}(0) \neq \emptyset$, and that   the   cocoercive parameter
$\lambda$ and the damping parameter $\gamma$ satisfy
 \be
  \lambda{\gam}^2 > 1.
 \ee
  Then, for each  initial data $u_0$ and $v_0$ in $\Hs$,
the unique  solution $u \in C^2([0,+\li); \Hs)$ of (\ref{pbn})
satisfies:

i1)  There exists $u_{\infty}\in  S $  such that  $u(t)
\rightharpoonup u_{\infty}$ weakly in $\Hs$ as $t \to +\li$.

\noindent Moreover,

i2) $\udp \in L^2 (0, +\infty; \Hs)$; \  $\lim_{t \to +\infty} |\udp
(t)| = 0$;

i3) $\uds+\nabla \phi(u)+ Ap  \in L^{2} (0, +\infty; \Hs)$ whenever
$p \in S$;

i4) for every $p \in S$, \  $ \lim_{t \to + \infty} |u(t) - p |$ \
exists.

\end{theorem}
\ind Before  proceeding  with the proof of Theorem \ref{wac}, we
 establish  three  technical lemmas.
The following lemma plays a crucial role  in the proof of
optimality of the weak cluster points of a given trajectory of
(\ref{pb0}). Because of its independent interest, it is stated in a  general setting, the potential $\Phi$ being allowed to be non-smooth.
We use the classical notation $\der \Phi$ for the subdifferential operator of $\Phi$ (it coincides with the gradient operator in the smooth case).
\begin{lemma} \label{opt}
 Let $\Phi : \Hs \to \R \cup \left\{+\infty\right\}$ be a convex and  lower semicontinuous
 function,
    $(u_n)$  a bounded  sequence in $\Hs$, and $p \in \Hs$.  If  there exist $\eta_n \in \der \Phi(u_n)$  and
     $\bar{\eta}\in \der \Phi(p)$  such that
 \be \label{nur1} \lim_{n \to +\li} \la u_n - p, \eta_n - \bar{\eta} \ra
 =0, \ee
 then  any weak cluster point $\bar{u}$ of $(u_n)$ satisfies
 \be \label{nur2}   \bar{\eta} \in \der \Phi(\bar{u}) . \ee
 \end{lemma}
 \ddem
Set $w_n= \la u_n - p, \eta_n - \bar{\eta} \ra $ and introduce the
functional $F:\Hs \to \R \cup \left\{+\infty\right\}$ defined for any $v \in \Hs$ by
 \be  \label{fff} F(v)= \Phi(v)- \Phi(p) -\la v-p, \bar{\eta}  \ra.
 \ee
 By   convexity of $\Phi$ and $\bar{\eta}\in \der
\Phi(p)$, we observe that $F$ is a convex and nonnegative function,
and it can be easily checked
 that
  \[  F(u_n)= \Phi(u_n)-\Phi(p) - \la u_n-p, \eta _n  \ra +  w_n. \]
Therefore, from  convexity of $\Phi$ and $\eta_n \in \der
\Phi(u_n)$, we have
\[  0 \le  F(u_n) \le   w_n, \]
 which  by  (\ref{nur1}) (that is $w_n
\to 0$  as $n \to +\li$) gives
 \be \label{nur3} \lim_{n \to +\li}
F(u_n)= 0 . \ee
 Now   consider  any   weak cluster point $\ov{u}$
of $(u_n)$, namely there exists a subsequence $(u_{n_k})$ of ($u_n)$
such that $(u_{n_k})$ weakly converges to $\ov{u}$  as $k \to +\li$.
Then invoking the weak lower semicontinuity of $F$ (as it is
convex and lower semicontinuous) and using (\ref{nur3}), we obtain

\begin{center}
$ 0 \le F(\ov{u}) \le \liminf_{k  \to +\li} F(u_{n_k}) = \lim_{n \to +\li} F(u_n)= 0,$
\end{center}

\noindent which entails  $ F(\ov{u})=0.$ This implies
that $\ov{u}$ is a minimizer of the convex and nonnegative function
$F$, so  that $ 0 \in \der F(\ov{u})$,
  which from (\ref{fff}) is equivalent to
(\ref{nur2})\edem

\begin{lemma} \label{wcp}
Let $A: \Hs \to \Hs$ be  a maximal monotone single-valued operator
and  $\phi: \Hs \to \R$ a convex differentiable function such that
$S:=( \nabla \phi +A)^{-1}(0) \neq \emptyset$. Suppose that  $p \in
S$ and $(u_n)$ is  a bounded sequence
 in $\Hs$ verifying
\[ \ba
{\rm (c1)} \es \lim_{n \to \li}  \la \nabla \phi(u_n)-\nabla \phi(p), u_n -p \ra =0, \\
{\rm (c2)} \es \lim_{n \to \li}  |Au_n -Ap| =0. \ea
 \]

\noindent Then, any weak cluster point of $(u_n)$ belongs to $S$.
\end{lemma}
\ddem Let  $\bar{u}$ be a  weak cluster point of $(u_n)$. From (c1)
and invoking Lemma \ref{opt}, we have $\nabla \phi(p)=\nabla
\phi(\bar{u})$. Moreover, from (c2)  and recalling that the maximal
monotone operator $A$ has a graph which is closed in $w-\Hs \times
s-\Hs$ (see \cite{Bre} for example), we have $A\bar{u}=Ap$. As a
straightforward consequence,  $\nabla \phi(\bar{u}) + A \bar{u}=
\nabla \phi(p) + Ap = 0$, so that
 $\bar{u} \in S$, which completes the proof\edem

\smallskip

The following  lemma  will be  used in the proof of Theorem
\ref{wac} in order to obtain the asymptotic convergence of the mapping $t
\to |u(t)-p|$, whenever $p \in S$ and $u$ is  solution of
(\ref{pb0}). This lemma appears implicitly in \cite{alv}, its
proof is given for the sake of completeness.
\begin{lemma} \label{ralv}
If $w \in C^2([0,+\li);\R)$ is bounded from below and  satisfies the
following inequality
\be \label{ral} \ddot{w}(t)+ \gam \dot{w}(t) \le g(t), \ee
where  $\gamma$ is a positive constant  and $g \in L^1
([0,+\li);\R)$, then $w(t)$ converges as $t \to + \li$.
\end{lemma}
\ddem From  (\ref{ral}) and denoting $[\dot{w}]_+=\max\{ \dot{w}(t), \ 0 \}$, we classically  have

\begin{center}
$ [\dot{w}(t)]_+ \le e^{-\gam t} [\dot{w}(0)]_+  + \int_0^t  e^{-\gam( t-\tau)} |g(\tau)| d\tau$,
\end{center}

\noindent while  Fubini's theorem gives us

 $$ \int_0^{+\li} \int_0^t e^{-\gam( t-\tau)}|g(\tau)| d \tau dt= \frac{1}{\gam} \int_0^{+\li}|g(\tau)| d\tau< +\li. $$

\noindent This  shows that $[\dot{w}]_+ \in L^1 ([0,+\li);\R)$. Now
setting $ z(t) = w(t) - \int_0^t [\dot{w}(\tau)]_+ d \tau$,  we  observe
that $z(.)$ is bounded from below (as $w(.)$ is assumed to be
bounded from below) with  $\dot{z}(t) = \dot{w}(t) - [\dot{w}(t)]_+ \le 0$,
hence $z(t)$ converges as $t \to +\li$, and so does $w(t)$\edem

\smallskip

\ind We are now in position to prove the main result of this
section.

\smallskip

\noindent {\bf Proof of Theorem \ref{wac}}.
  Take $p\in S = ( \nabla
\phi + A)^{-1}(0)$ and set $h(t)=(1/2) |u(t)-p|^2$. From $ \dot{h}(t) =
\la u(t)-p, \udp (t) \ra$ and $\ddot{h}(t) \dis = \la u(t)-p, \uds
(t) \ra + |\udp (t)|^2$  we obtain
  \be \ddot{h}(t) +  \gamma
\dot{h}(t) =  \la u(t)-p, \uds (t) + \gamma \udp (t) \ra  +  |\udp
(t)|^2, \ee
 which, by using (\ref{pbn}) yields
\be \ddot{h}(t) +  \gamma \dot{h}(t) +  \la  \nabla  \phi
(u(t)) + Au(t) , u(t)-p \ra = |\udp (t)|^2. \ee
Recalling that $ \nabla  \phi(p) + Ap = 0$, we deduce that
\be \label{supl4}  \ddot{h}(t) +  \gamma \dot{h}(t) + \la  \nabla  \phi (u(t)) - \nabla  \phi(p) , u(t)-p \ra   +  \la Au(t) -Ap, u(t)-p \ra = |\udp (t)|^2. \ee
 Then, by  using the
$\lambda$-cocoercive property (\ref{cc}) of $A$,  together with the
monotonicity of $\nabla  \phi$,  we obtain
\be \label {bsinq1} \ddot{h}(t) +  \gamma \dot{h}(t) +  \lambda
|Au(t)-Ap|^2  \leq |\udp (t)|^2. \ee
Regarding the third term in the left side of (\ref{bsinq1}), by
using (\ref{pbn}) again,  replacing $Au(t)$ by $-( \uds (t) +\gam
\udp(t) + \nabla \phi (u(t)))$, and setting
  $ D(t):=|\uds (t)+ \nabla \phi (u(t))+Ap|^2$,
we obtain
\begin{eqnarray*}
|Au(t)-Ap|^2 & = & |\lb \uds (t)
 + \nabla  \phi (u(t))+Ap \rb + \gam \udp(t)|^2\\
 & = & D(t) + \gam^2 |  \udp(t)|^2+  2 \gam \la \udp(t), \uds
(t)
 + \nabla  \phi (u(t))+Ap \ra \\
  & = & D(t)  + \gam^2 |  \udp(t)|^2 +  \gam
\frac{d}{dt}\lb |\udp(t)|^2 + 2 {\phi}(u(t)) + 2 \la u(t)-p, Ap \ra \rb.
\end{eqnarray*}

\noindent Recalling that $Ap=-\nabla \phi(p)$ and by using the above equality,
let us   rewrite (\ref{bsinq1}) as

 \be \label {bsinq2} \ba
 \dis  \ddot{h}(t) +  \gamma
\dot{h}(t)+ (\lambda {\gamma}^2 - 1) |\udp (t) |^2
 + \lambda D(t) \\
\hskip 2.cm \dis   +  \lambda \gamma \frac{d}{dt} \lb |\udp(t)|^2 +
2 {\phi} ( u(t)) - 2 \la u(t)-p, \nabla \phi(p) \ra \rb  \leq 0. \ea
\ee
By using assumption $\lambda {\gamma}^2 \geq 1$, with $D(t) \ge 0$
and (\ref{bsinq2}), we then get
\be \label {bsinq3} \ddot{h}(t) +  \gamma \dot{h}(t) +  \lambda
\gamma \frac{d}{dt}\lb |\udp(t)|^2 + 2 {\phi}(u(t))  - 2 \la u(t)-p,  \nabla
\phi(p) \ra \rb \leq 0, \ee
 which expresses that the function
\be \label {bsinq4} \Gamma_0 (t) := \dot{h}(t) +  \gamma h(t) +
\lambda \gamma \lb |\udp(t)|^2 + 2 {\phi}(u(t)) - 2 \la u(t)-p, \nabla
\phi(p) \ra \rb \ee
\noindent is nonincreasing on $\left[0, +\infty \right)$. Indeed,
$\Gamma_0(\cdot) $ will serve us as a Liapunov function in  the
asymptotic analysis of (\ref{pbn}).  Let us first show the
boundedness of  trajectories. Note that by   convexity  of $\phi$,
we obviously have
\be \label{ou1}  {\phi}(p) \le  {\phi}(u(t)) -  \la u(t)-p, \nabla
\phi(p) \ra, \ee
hence  the nonincreasing property of $\Gamma_0(\cdot)$ leads to
\be \label {bsinq31}
 \dot{h}(t) +  \gamma h(t) +
\lambda \gamma \lb  |\udp(t)|^2 +  2 {\phi}(p) \rb \le
\Gamma_0(t)\le \Gamma_0(0).
\ee
Set

\begin{center}
$C_0:= |u_0 - p||v_0| + \frac{\gamma}{2}|u_0 - p|^2 + \lambda \gamma
|v_0|^2 +  2 \lambda \gamma \left({\phi}(u_0)-\phi(p)- \la u_0-p,
\nabla \phi(p) \ra \right) $,
\end{center}

\noindent which clearly satisfies

\be \label{bins1} C_0 \geq \Gamma_0(0)-2\lambda \gamma \phi(p). \ee

\noindent From (\ref{bsinq31}) and (\ref{bins1}) we
 deduce that
%
$ \dot{h}(t) +  \gamma  h(t)   \leq C_0, $
%
so that,  after  integration, we obtain

\begin{center}
$ h(t)= \frac{1}{2}|u(t) - p|^2 \leq \frac{1}{2}|u_0 - p|^2 +  \frac{C_0}{\gamma},$
\end{center}

\noindent which  shows that the trajectory $t \to u(t)$ remains bounded on
$\left[0, +\infty \right)$:

\be \label {bsinq6c}
 \sup_{t\in \left[0, +\infty \right)} |u(t)| < + \infty .
\ee

 Let us now establish  estimates  on $\udp$ and $\uds$. Integrating
  (\ref{bsinq2})  from 0 to $t$,
together with (\ref{ou1}) yields
 \be \label{supl2} \ba
 \dis  \dot{h}(t) +  \gamma h(t) +  \lambda \gamma |\udp(t)|^2 \\
 \hskip 1.cm \dis + (\lambda {\gamma}^2 - 1) \int_0^t |\udp (s)|^2
ds + \lambda \int_0^t |\uds (s)+ \nabla \phi(u(s))+Ap|^2 ds \leq
C_0, \ea \ee
which implies
%
$\dot{h}(t) +   \lambda \gamma |\udp(t)|^2  \leq C_0$.
%
Equivalently, by definition of $h(t)$, we have

   \be \label {bsinq8} \la
u(t)-p, \udp(t) \ra   + \lambda \gamma |\udp(t)|^2  \leq C_0. \ee
\noindent Since $u(t)$ remains bounded on $\left[0, +\infty
\right)$, inequality (\ref {bsinq8})
 implies that $\udp(t)$ also remains bounded on $\left[0, +\infty \right)$:
\be \label {bsinq9}
\sup_{t\in \left[0, +\infty \right)} |\udp(t)| < + \infty .
\ee
From (\ref{bsinq6c}),  (\ref{bsinq9}) and $\dot{h}(t) = \la u(t)-p,
\udp (t) \ra$, we observe  that $\dot{h}(t)$ is bounded on $\left[0,
+\infty \right)$. Returning to (\ref{supl2}) we deduce that
\be \label {bsinq10}  \label{bsinq11} \mbox{$\dis \int_0^{+ \infty}
|\udp (s)|^2 ds < + \infty$ \ and \ $\dis \int_0^{+ \infty} |\uds (s)+
\nabla \phi (u(s)) +Ap|^2 ds < +\infty$.} \ee
\noindent Using that  $u$ and  $\udp$ are bounded on  $\left[0,
+\infty \right)$, together with the Lipschitz continuity property of $A$  and equation (\ref{pb00}), we
deduce that $\uds$ is also bounded:
\be \label {bsinq12}
\sup_{t\in \left[0, +\infty \right)} |\uds(t)| < + \infty .
\ee
\noindent Properties (\ref{bsinq10}) and (\ref{bsinq12}) express that the function $g:= \udp$ satisfies both
\begin{center}
$g \in L^2 (0, +\infty; \Hs)$  \ and \ $ \dot{g} \in L^{\infty} (0, +\infty; \Hs)$.
\end{center}
\noindent By a classical result, this implies  $\lim_{t\to  +\infty} g(t) = 0 $, that is,
\be \label {bsinq13} \lim_{t\to  +\infty}| \udp(t)| = 0 . \ee
\noindent Returning to (\ref{pb00}), by using that functions  $u$
and $\udp$ are Lipschitz continuous on $\left[0, +\infty \right)$
(as their derivatives are uniformly bounded), that $A$ is Lipschitz
continuous and $\nabla\phi$ is Lipschitz continuous on bounded sets,
we deduce that  $\uds$ is  Lipschitz continuous on $\left[0, +\infty
\right)$. Once again, this property together with  (\ref{bsinq10})
implies that
$ \label {bsinq14} \lim_{t\to  +\infty} | \uds(t)+ \nabla \phi(u(t))
+ Ap |  = 0$,
 which  in light of  (\ref{pb00}) and (\ref{bsinq13}) leads to
 \be \label {supl3} \lim_{t\to  +\infty} |A u(t)-Ap|=0.  \ee
\ind We have all the ingredients  to conclude thanks to Opial's
lemma, see \cite{Op}. To that end, we  need to prove  that
properties (a) and (b) are fulfilled:

(a) for every $p \in S$, \  $ \lim_{t \to + \infty} |u(t) - p |$ \
exists;

(b) for every $t_n \to + \infty$ with $ u(t_n) \rightharpoonup
\overline{u}$ weakly in $\Hs$, \  we have $\overline{u} \in S$.

\noindent Opial's lemma asserts that, under the  above two
properties,  $u(t)$ weakly converges as $t\to + \infty$ to an
element $u_{\infty} \in S$.

\noindent Let us first prove (a).
  From
(\ref{supl4}) and the monotonicity of $A$,  we  have
 \be \label{supl5} \ddot{h}(t)
+ \gamma \dot{h}(t) + w(t) \le  |\udp (t)|^2 , \ee
where \[  w(t)= \la \nabla \phi (u(t)) - \nabla \phi(p) , u(t)-p \ra
. \] It is obviously seen that $w(t)$ is   a nonnegative term
(thanks to the monotonicity of $\nabla \phi$), hence,
  from (\ref{supl5}), we immediately derive
 \be \ddot{h}(t)
+ \gamma \dot{h}(t)  \le  |\udp (t)|^2 . \ee
By (\ref{bsinq10}) we know that $|\udp |^2$ belongs to $L^1(0,+\li;
\Hs)$. Noticing that  $h$ is nonnegative,
Lemma \ref{ralv} shows  that property (a) holds.

\noindent  Let us now prove (b). By (\ref{supl5}),  after integration we obtain
 $\int_0^{+\li} w(s) ds < + \li$. Owing to the Lipschitz continuity property of $w(\cdot)$  we deduce that
  $ \lim_{ t \to + \li} w(t) = 0$. On the other hand,  by (\ref{supl3})  we
 have  $\lim_{t \to +\li} |Au(t)-Ap| = 0$. We can now apply  Lemma \ref{wcp} to get property (b)\edem
\medskip

\subsection{The strongly monotone case}

 Let us recall that  the operator $\nabla \phi + A $ is called {\it strongly monotone over bounded
 sets} if,  for all positive real number $R$, there exists a
 continuous function $w_R: \R_+ \to \R_+$ verifying
 \[
   w_R(t_n)\to 0 \Rightarrow t_n \to 0, \]
  and   such that  $(u,v) \in \Hs^2$ with $|u| \le R$ and $|v| \le
  R$ yields
 \be \label{smbs}
 \la  (\nabla \phi + A)u - ( \nabla \phi + A)v, u -v \ra \ge w_
 {R}(|u-v|).
 \ee

 \noindent Clearly,
  this property implies that  $S:=( \nabla \phi + A)^{-1}(0)$ contains at most one element and it
  holds  if, for instance,   $\nabla \phi \ $ or $A$ is strongly monotone over bounded
  sets (with the other operator being assumed to be monotone).

\begin{proposition}\label{prop1}
Under the  assumptions of  Theorem \ref{wac}, assuming moreover that  $ \nabla \phi +A$
 is  strongly monotone over bounded sets, then
the solution $u(t)$ of (\ref{pb0})  strongly converges as $t$ goes to $+\infty$  towards
 the unique element  of  $\ S:=( \nabla \phi + A)^{-1}(0)$.
\end{proposition}
\ddem By Theorem \ref{wac}, we have $ u, \dot{u}  \in
L^{\li}([0,+\li);\Hs)$. Set
 $R:= \sup_{ t \ge 0} |u(t)|+|p|$, with $p \in S$, and combine  (\ref{supl4}) and (\ref{smbs}),
 to obtain
  $ \ddot{h}(t)  +  w_R(|u(t)-p|) \le  |\udp (t)|^2$,
   where $h(t):=(1/2)|u(t)-p|^2$. After integration, observing that
   $\dot{h} \in L^{\li}([0,+\li);\Hs)$ and recalling that $\dot{u}
   \in L^{2}([0,+\li);\Hs)$ (Theorem \ref{wac}, (i2)), we deduce
   that $\int_0^{+\li}w_R(|u(s)-p|)ds < +\li$. On the other hand, by Theorem \ref{wac}, (i4), $\lim_{ t \to +\li}|u(t)-p|$ exists, which by continuity of $w_R$ implies that
   $\lim_{ t \to +\li}w_R(|u(t)-p|)$ exists.  Hence,     $\lim_{ t \to +\li}
   w_R(|u(t)-p|) = 0$.  Then use the property of $w_R$ ($ w_R(t_n)\to 0 \Rightarrow t_n \to 0$) to obtain $\lim_{ t \to +\li}
   |u(t)-p|=0$, which ends the proof\edem


\subsection{Some particular cases}

Let us specialize Theorem  \ref{wac} and Proposition \ref{prop1} to
some important particular situations.

a) Let us first take  the non-potential operator $A$ equal to zero.
Then notice that the null operator is $\lambda$-cocoervive for any
$\lambda>0$. By taking $\lambda > (1/\gamma^2)$, we can apply
Theorem  \ref{wac}. Owing to the fact that no restrictive assumption
is made on the potential part, we obtain  the Alvarez theorem
\cite{alv}:
\begin{corollary} \label{cor1}    Let  $\phi : \Hs \to \R$ be a
 convex differentiable function whose gradient $\nabla \phi$ is
 Lipschitz continous on the bounded subsets of $\Hs$, and let $\gamma >0$.
 Let us assume that $S:= \argmin_{\Hs}\phi$ is nonempty.
Then, for each  initial data $u_0$ and $v_0$ in $\Hs$,
the unique  solution $u \in C^2([0,+\li); \Hs)$ of

\be  \label{cor11}
  \uds (t) +\gam \udp(t)   + \nabla \phi(u(t))  = 0,
 \ee

\noindent with initial conditions $u(0)=u_0$ and  $\dot{u}(0)=v_0$,
satisfies:

i1)  there exists $u_{\infty}\in  S$  such that  $u(t) \rightharpoonup u_{\infty}$ weakly in $\Hs$ as $t \to +\li$;

i2) $\udp \in L^2 (0, +\infty; \Hs)$; \  $\lim_{t \to +\infty} |\udp
(t)| = 0$.

\end{corollary}

b) Operators of the form  $A = I - T$, where $T: \Hs \to  \Hs$ is a
contraction, play a central role in the realm of fixed point theory
and constrained optimization. Let us verify  that $A$ is
$(1/2)$-cocoercive.
 Indeed, by setting
\bc
  $ E := \left\langle  (u -Tu) - ( v-Tv), u-v \right\rangle   - \frac{1}{2}|(u -Tu) - ( v-Tv)|^2  $,
\ec
  we obviously have

\[ \ba E=  \left\langle  (u -v) - ( Tu-Tv), u-v \right\rangle - \frac{1}{2}|(u -v) - ( Tu-Tv)|^2
 , \\
 \mbox{ } \mbox{ }=  |u -v|^2  - \left\langle  Tu - Tv, u-v \right\rangle  - \frac{1}{2}|u -v |^2  - \frac{1}{2}| Tu-Tv|^2 + \left\langle  (u -v) ,  Tu-Tv \right\rangle, \\
 \mbox{ } \mbox{ } =  \frac{1}{2}(|u -v |^2  - | Tu-Tv|^2),
 \ea
 \]
which, by  contraction property of $T$, is nonnegative. Thus we can
take $ \lambda = 1/2$, and  condition (\ref{cf}) boils down to
$\gamma
> \sqrt{2}$. Applying Theorem  \ref{wac}  to this situation, we
obtain the following result (see \cite {alat2}, Theorem 3.2):

\begin{corollary} \label{cor2}    Let  $T: \Hs \to  \Hs$ be a contraction and  $\gamma > \sqrt{2}$.  Let us assume that
 $S:= Fix T = \left\{v\in \Hs: Tv = v \right\}$ is nonempty.
Then, for each  initial data $u_0$ and $v_0$ in $\Hs$,
the unique  solution $u \in C^2([0,+\li); \Hs)$ of

\be  \label{cor22}
  \uds (t) +\gam \udp(t)   + u(t) - T(u(t))  = 0,
 \ee

\noindent with initial data $ u(0)=u_0$ and  $\dot{u}(0)=v_0$,
satisfies:

i1)  there exists $u_{\infty}\in  Fix T$  such that  $u(t) \rightharpoonup u_{\infty}$ weakly in $\Hs$ as $t \to +\li$;

i2) $\udp \in L^2 (0, +\infty; \Hs)$; \  $\lim_{t \to +\infty}
|\udp(t)| = \lim_{t \to +\infty} |\uds(t)| = 0$.

\end{corollary}

c) When working with a general maximal monotone operator $A: \Hs \to  \Hs$, let us observe that, for every $\lambda >0$, its Yosida approximation
 $A_{\lambda} $ is $\lambda$-cocoercive.
We recall that   $ A_{ \lambda }= \frac{1}{\lambda}(I -J_{\lambda}^A)$, where $J_{\lambda}^A:= (I +
\lambda A)^{-1}$
 is the resolvent operator of index $\lambda$ of $A$. The operator  $J_{\lambda}^A$  is everywhere  defined,  single-valued (it is a contraction) and
 $ A_{ \lambda} v  \in (A \circ J_{\lambda}^A)v$ for any $v \in \Hs$ (see \cite{Bre}).
 Let $(u,v) \in \Hs^2$ and set $E:= \left\langle A_{ \lambda }u - A_{ \lambda }v , u-v \right\rangle$.
  Noticing that
 $u= J_{\lambda}^A u + \lambda A_{ \lambda}u$ and  $v= J_{\lambda}^A v + \lambda A_{ \lambda}v$,
 and recalling that
 $ A_{ \lambda} u  \in A  (J_{\lambda}^Au)$ and $ A_{ \lambda} v  \in A  (J_{\lambda}^Av)$,
 we immediately obtain

 \[   E =  \left\langle   A_{ \lambda }u - A_{ \lambda }v , (J_{\lambda}^A u -  J_{\lambda}^A v)
 + \lambda (A_{ \lambda}u - A_{ \lambda}v) \right\rangle
 \mbox{ } \mbox{ }\geq  \ \lambda |A_{ \lambda }u - A_{ \lambda }v
 |^2,
 \]
 which proves that $A_{\lambda} $ is $\lambda$-cocoercive.
  Observing  that $A$ and
$A_{\lambda}$ have the same zeroes, and as a straight consequence of
Theorem \ref{wac}, we reach  the following result:

\begin{corollary} \label{cor3}    Let  $A: \Hs \to  \Hs$ be a
general maximal monotone operator and let $\gamma >0$ and $\lambda>0$ be such that $\lambda{\gam}^2 > 1$.
 Let us assume that  $S:= A^{-1} (0)$, the set or zeroes of  $A$, is nonempty.
Then, for each  initial data $u_0$ and $v_0$ in $\Hs$,
the unique  solution $u \in C^2([0,+\li); \Hs)$ of

\be  \label{cor33}
  \uds (t) +\gam \udp(t)   + A_{\lambda}(u(t))  = 0,
 \ee

\noindent with initial data  $u(0)=u_0$ and $ \dot{u}(0)=v_0$,
$A_{\lambda}$ being  the Yosida approximation of index $\lambda$ of
$A$, satisfies:

i1)  there exists $u_{\infty}\in  S$  such that  $u(t) \rightharpoonup u_{\infty}$ weakly in $\Hs$ as $t \to +\li$;

i2) $\udp \in L^2 (0, +\infty; \Hs)$; \  $\lim_{t \to +\infty} |\udp
(t)| = 0$.

\end{corollary}

Let us examine an interesting consequence of Corollary \ref{cor3}
regarding numerical schemes. From  resolvent equation  (see
\cite{Bre}, Proposition 2.6), we have  $ (A_{\lambda})_{\mu} =
A_{\lambda  + \mu}$ whenever  $\lambda
>0$ and  $\mu > 0$.  Then it  can be easily derived that,  for any $v\in \Hs$,

\begin{equation}\label{resol}
J_{\mu}^{A_{\lambda}}v = \frac{\lambda }{ \lambda + \mu}v  +  \frac{\mu}{ \lambda + \mu}  J_{\lambda + \mu}^{A}v .
\end{equation}

\noindent Implicit time discretization of (\ref{cor33})   naturally leads to a second order relaxed proximal algorithm, whose trajectories (sequences)
 converge to the set of zeroes $A^{-1} (0)$, for $A$ a general maximal monotone operator, see \cite{alv2} and references herein.

\subsection{Condition (\ref{cf}) is sharp}

 Let $B: \R^2  \to \R^2$ be  the $\pi/2$ rotation of center
$(0,0)$, namely  $B$ is the linear operator whose
matrix in the canonical basis is given by
\[ B=\left ( \begin{array}{cc}
0 & -1 \\
1 &  0 \\
\end{array} \right ).
  \]
Consider the  dynamical system
 \be \label{sexp1} \ddot{X}(t) + \gamma \dot{X}(t) +B_{\ld}X(t)=0 ,
 \ee
where $\gamma>0$ and $B_{\ld}$  is the Yosida approximation
 of $B$ of parameter $\lambda  >0$. An easy computation shows that

\[ B_{\ld} = \frac{1}{1+ \ld^2} H , \es \mbox{where} \es H= \left ( \begin{array}{cc}
\ld  & -1 \\
1 &  \ld \\
\end{array} \right ).
  \]

\noindent The operator $B_{\ld}$ is $\lambda$-cocoercive, as follows from general properties of Yosida approximation, or  by a direct elementary computation.
Let us compute explicitly the solutions of (\ref{sexp1}).
The matrix $H$ admits two complex conjugate eigenvalues. One
  of them is
  $\rho_1=\ld-i$, and it is associated to  the     eigenvector
  $ W_1=(1,i)^T$.
Subsequently, searching for a solution of (\ref{sexp1}) of the form
$X(t)= e^{rt} W_1$ (with $r \in \C$), we see that $r$ satisfies the
quadratic equation
 \be \label{esd} r^2 + \gamma r + \frac{\ld-i}{1+ \ld^2}=0. \ee
Equation (\ref{esd})  has a complex
discriminant $  \Delta= \gamma^2 - \frac{4 \ld
}{1+\ld^2}+\frac{4i}{1+\ld^2}  =(x+iy ) ^2$,
 where $x$ and $y$ are real values given by
 \[ \ba
 x = \frac{1}{\sqrt{2}} \lb (\gamma ^2 - \frac{4 \ld}{1+ \ld^2} )
 + \sqrt{ \lb \gamma^2-\frac{4 \ld }{1+\ld^2} \rb ^2
 + \frac{16}{(1+ \ld^2)^2}} \rb^{1/2}, \\
 y = \frac{1}{\sqrt{2}} \lb -(\gamma ^2 - \frac{4 \ld}{1+ \ld^2} )
 + \sqrt{ \lb \gamma^2-\frac{4 \ld }{1+\ld^2} \rb ^2
 + \frac{16}{(1+ \ld^2)^2}}
 \rb^{1/2}.
  \ea
 \]
 Hence (\ref{esd}) has  two complex
solutions
 \[ \ba r_1=(1/2)[(-\gamma -x)-iy]= a_1-ib, \\
r_2=(1/2)[(-\gamma+x )+ iy]= a_2+ib, \ea \]
 where $a_1$, $a_2$ and
$b$ are real values defined by
 \[ \ba a_1=\frac{1}{2} \lb -\gamma -
  \frac{1}{\sqrt{2}} \lb (\gamma ^2 - \frac{4 \ld }{1+ \ld^2} )
 +\sqrt{ \lb \gamma^2-\frac{4 \ld }{1+\ld^2} \rb ^2
 + \frac{16}{(1+ \ld^2)^2}} \rb^{1/2}
 \rb, \\
 a_2=\frac{1}{2} \lb -\gamma + \frac{1}{\sqrt{2}} \lb (\gamma ^2 - \frac{4 \ld }{1+ \ld^2} )
 + \sqrt{ \lb \gamma^2-\frac{4 \ld }{1+\ld^2} \rb ^2
 + \frac{16}{(1+ \ld^2)^2}} \rb^{1/2}
 \rb, \\
b=  \frac{1}{2\sqrt{2}} \lb (-\gamma ^2 +\frac{4 \ld}{1+ \ld^2} )
 +\sqrt{ \lb \gamma^2-\frac{4 \ld }{1+\ld^2} \rb ^2
 + \frac{16}{(1+ \ld^2)^2}}
 \rb^{1/2}.
 \ea
 \]
  As a
  straightforward consequence, we deduce that
a family  of complex solutions to (\ref{sexp1}) is given by the
functions $\{ V_1,
 V_2\}$, where $V_1(.)$ and $V_2(.)$ are
 defined for any $ t \in \R$ by
 \[ \ba  V_1(t)=
 \left ( \begin{array}{c}
1 \\
i \\
\end{array} \right ) e^{(a_1-ib)t}  \\
= e^{a_1 t}
\left ( \begin{array}{c}
\cos bt - i \sin bt  \\
 i \cos bt + \sin bt \\
\end{array} \right )
= e^{a_1 t} \lp \left ( \begin{array}{c}
\cos bt  \\
  \sin bt   \\
\end{array} \right ) + i \left ( \begin{array}{c}
 -\sin bt  \\
  \cos bt  \\
\end{array} \right ) \rp,
\ea  \]
\[ \ba  V_2(t)=
 \left ( \begin{array}{c}
1 \\
i \\
\end{array} \right ) e^{(a_2+ ib)t}  \\
= e^{a_2 t} \left ( \begin{array}{c}
\cos bt + i \sin bt  \\
 i \cos bt - \sin bt \\
\end{array} \right )
= e^{a_2 t} \lp \left ( \begin{array}{c}
\cos bt  \\
 - \sin bt   \\
\end{array} \right ) + i \left ( \begin{array}{c}
 \sin bt  \\
  \cos bt  \\
\end{array} \right ) \rp.
\ea  \]
 Note also that the complex conjugate of any solution to
  (\ref{sexp1}) is also a solution to  (\ref{sexp1}).
 It is then immediate that a family  of real solutions to (\ref{sexp1}) is
given by the  functions  $\{ U_1,
 U_2,  U_3, U_4 \}$ defined for any $ t \in \R$ by
 \[  U_1(t)=
e^{a_1 t} \left ( \begin{array}{c}
 \cos bt \\
 \sin bt  \\
\end{array} \right ) , \es  U_2(t)=
e^{a_1 t} \left ( \begin{array}{c}
- \sin bt \\
 \cos bt  \\
\end{array} \right ), \]
\[  U_3(t)=
e^{a_2 t} \left ( \begin{array}{c}
 \cos bt \\
 -\sin bt  \\
\end{array} \right ) , \es  U_4(t)=
e^{a_2 t} \left ( \begin{array}{c}
 \sin bt \\
 \cos bt  \\
\end{array} \right ). \]
In light of these last results, we deduce that (\ref{sexp1})
will have non-convergent solutions  if  $a_2 \ge 0$, or equivalently
if
  \be \label{cir1}    \sqrt{ \gamma^4
 -  \frac{8 \gamma^2 \ld}{1+ \ld^2} + \frac{16}{1+ \ld^2}}
  \ge    (\gamma ^2 +
\frac{4 \lambda}{1+ \ld^2} ).
        \ee
 In addition, setting $ \ld = \theta / \gamma^2$ (with   $\theta \ge 0$), and after an easy
 computation, we obtain that (\ref{cir1}) is equivalent to
 $ \gamma^4 (1-\theta) \ge \theta^3$.
Hence (\ref{sexp1}) has   non-convergent
  solution trajectories  when $\theta \in [0,1)$, which means that
  $\ld \gamma^2 < 1$ does not ensure  convergence of  dynamical system
  (\ref{pb0}).

\subsection{Invariance of condition (\ref{cf})  with respect to time rescaling}

Given some positive real parameter $k$, let us consider the  time
rescaling $t=ks$. For any trajectory $u(.)$ of (\ref{pb0}), the
rescaled trajectory $v(s) = u(ks)$ satisfies
 \be \label {rescale1}
\ddot{v}(s) + \gam k\dot{v}(s) + k^2\nabla \phi(v(s)) +  k^2 A(v(s))
= 0. \ee
\noindent Let us suppose that $A$ is a $\lambda$-cocoercive maximal monotone operator.
 One can easily verify that $k^2 A$ is a ($\lambda / k^2)$-cocoercive maximal monotone operator.
\noindent Clearly, this time rescaling does not change the asymptotic convergence properties of the
trajectories. Indeed, the condition insuring the convergence of the trajectories,
namely  $\lambda{\gam}^2 > 1$, remains invariant under this time rescaling, as shown by
the following relation:
%
$(\lambda / k^2)  \times (\gam k)^2 = \lambda{\gam}^2.$
%
\noindent This elementary observation tells us that the condition
$\lambda{\gam}^2 > 1$ makes sense from a physical point of view.

 \section{\large Asymptotic stabilization by Tikhonov regularization methods}
 In order to correct some drawbacks of system (\ref{pb00}), namely the weak (and possibly not strong) asymptotic convergence property of its trajectories,
 and the dependence of the limit equilibrium on the initial data, we introduce in the equation a Tikhonov-like regularization term with a vanishing coefficient. To be more precise,
we are going to establish strong asymptotic convergence results
 regarding the solution of the non-autonomous system
 \be \label{pb1}
 \left [ \ba
  \uds (t) +\gam \udp(t)  + \nabla \phi(u(t))  +   A(u(t)) + \eps(t) \nabla \Theta (u(t)) = 0,\\
  \\
 u(0)=u_0, \es \dot{u}(0)=v_0, \ea \right . \ee
where $\gamma >0$ and $(u_0,v_0)$
in $\Hs^2$  are arbitrary given initial data. We make the   following
assumptions:

 (H1)  \begin{minipage}[t]{12.cm} $\phi : \Hs \to \R$ is a
 convex differentiable function whose gradient $\nabla \phi$ is Lipschitz continous on bounded sets;
\end{minipage}

\smallskip

  (H2) \begin{minipage}[t]{12.cm}
 $A: \Hs \to \Hs$ is maximal
monotone and $\ld$-cocoercive  for some $\ld >0$;
\end{minipage}

\smallskip

  (H3)  $\eps: [0,+\li) \to (0,+\li)$ is of class  $C^1$ and tends to zero  as $t
\to +\li$;

\smallskip

  (H4) \begin{minipage}[t]{12.cm} $\Theta: {\cal H} \to \R $ is
 differentiable,   convex and bounded below on $\Hs$,  and  its derivative $\nabla \Theta$
 is   $\delta$-Lipschitz continuous (LC)
  and $\eta$-strongly monotone (SM), with $\delta > 0$ and $\eta >0$, i.e., \\
\hspace*{1.cm}  (LC) $|\nabla \Theta (u) -\nabla \Theta(v)| \le
\delta |u-v|$  \  $\forall
  u,v\in \Hs$, \\
 \hspace*{1.cm}  (SM) $\la \nabla \Theta (u) -\nabla \Theta(v), u-v \ra  \ge \eta |u-v|^2$ \
$\forall
  u,v \in \Hs$.
  \end{minipage}

\smallskip

\noindent At once, we claim the main result of this section.

\begin{theorem} \label{imp}
Let us suppose that   {\rm (H1)-(H4) } hold with  $S := (A+ \nabla
\phi)^{-1}(0) \neq \emptyset$ and with $  \ld \gamma^2  > 1$.
Let us assume moreover that $\eps (t)$ is decreasing and converges slowly to zero as $t
\to +\li$ in the following sense:

    \be
   \int_0^{+\li} \eps(s) ds = +\li.
    \ee

\noindent Then, for any  $u_0$ and $v_0$ in $\Hs$,
there exists a unique  solution $u \in C^2([0,+\li); \Hs)$ of (\ref{pb1})
which satisfies  $u(t) \to u_*$ strongly in  $\Hs$   as $t \to +\li$, where $u_*$ is the
unique solution of  the variational inequality: find $u_* \in S$ such that

\begin{center}
$ \la \nabla \Theta (u_*), v-u_* \ra \ge 0 \es \forall v\in S.$
\end{center}

 \end{theorem}
\ind Before proving  Theorem \ref{imp}, we introduce a series of
preliminary results.  To begin with, we establish a key estimate on
the trajectories of  (\ref{pb1}).
\begin{lemma} \label{res2}
 Under conditions  {\rm (H1)-(H4)},   the   solution $u$
 of (\ref{pb1}) satisfies
 \be \label{int1} \ba
 \ld G(t) + w(t)+ (\ld \gam ^2-1) | \udp(t)|^2 + \dot{\Gamma}_1(t) \le \\
  \hskip 2.cm  - \eps(t) \la \nabla {\Theta} (u(t)), u(t)-p \ra + 2
\ld \gam \dot{\eps}(t) \bar{\Theta}(u(t)), \ea  \ee
 where  $p$ is any element in $S:=(A+ \nabla \phi) ^{-1}(0)$, $w(.)$
 and $G(.)$ are defined by \\
  \ind $w(t)=\la  \nabla
 \phi(u(t))-\nabla \phi(p), u(t)-p \ra$,  $G(t):= | A(u(t)) + \gam \dot{u}(t) -Ap|^2$, \\
 $\dot{\Gamma}_1 (\cdot)$ denotes  the first derivative of
 the mapping $\Gamma_1(\cdot)$ given  by  \\
 \ind $ \Gamma_1(t)=\dot{h}(t) + \gam h(t) + \ld  \gamma  \lb  |\udp (t)|^2 +
 2   \lb  \phi(u(t)) - \la u(t)-p, \nabla \phi(p)\ra \rb  + 2  \eps(t)  \bar{\Theta}(u(t))
 \rb$, \\
 with   $h(t)=(1/2)|u(t)-p|^2$,
 $\bar{\Theta}(v):=\Theta(v)-\inf _{\Hs} \Theta$.
  \end{lemma}
   \ddem
   Without ambiguity, in order to get simplified notations, we omit to write the variable $t$. Hence $u$ stands for $u(t)$, and so on.
Taking  $p \in S$ and using (\ref{pb1}), we have
 \be \label{gis2}
   \la \nabla \phi(u)-\nabla\phi(p), u-p \ra  + \la Au-Ap, u-p \ra  = -\la  \uds + \gam \udp  + \eps(t) \nabla \Theta (u) ,
 u-p \ra,
 \ee
which by   $\ld$-cocoerciveness of $A$ and by definition of
$w(t)$ amounts to
  \be \label{tis5}
\ld   |Au-Ap  |^2 + w(t) \le -\la  \uds   + \gam \udp, u-p \ra -
\eps(t)  \la  \nabla \Theta (u), u-p \ra . \ee
 Furthermore, from $h=(1/2) |u-p|^2$, we have $
 \dot{h}=\la \udp, u-p \ra$ and $ \ddot{h}= \la \uds, u-p \ra+
 |\udp|^2$, so that
 $   \ddot{h} + \gamma \dot{h}= \la \uds + \gamma \udp , u-p \ra+   |\udp|^2
 $.
This combined with (\ref{tis5}) gives
 \be \label{bis5}
\ld   |Au-Ap  |^2 + w(t) + \lb \ddot{h} + \gamma \dot{h}  \rb \le
|\udp|^2-\eps(t) \la   \nabla \Theta (u),  u-p \ra. \ee

\noindent Let us reexpress $|Au-Ap|^2$ with the help of (\ref{pb1}):
 \be \label{bis51} \ba
 |Au -Ap |^2= | \uds  + \nabla \phi(u)+ \eps(t) \nabla \Theta (u) +Ap + \gam \udp |^2 \\
 \hskip 2.cm = G(t)+ \gam ^2 | \udp |^2 + 2 \gam \la \udp,
 \uds + \nabla \phi(u)+ \eps(t) \nabla \Theta (u)+Ap  \ra,
  \ea \ee

\noindent where  $ G(t):= | \uds + \nabla \phi(u)+ \eps(t) \nabla
\Theta (u)+Ap|^2$ (hence $G(t)= | Au(t)+ \gam \dot{u}(t) -Ap|^2$).

\noindent Set
$\bar{\Theta}(.)=\Theta(.) - \inf_{\Hs} \Theta$,
   $Q(t)= (1/2) |\udp|^2 + \phi(u) + \la u-p, Ap\ra  $,
and rewrite (\ref{bis51}) as
  \[  \ba \dis  |Au -Ap |^2  =  G(t)+ \gam ^2 | \udp |^2 + 2 \gam  \frac{d}{dt} \lb
   Q(t) + \eps(t)  \bar{\Theta}(u) \rb
   -2 \gamma  \dot{\eps}(t)  \bar{\Theta} (u). \ea
 \]
\noindent Replacing this last expression in  (\ref{bis5}) we obtain
 \[ \ba
w(t)+  (\ddot{h} + \gamma \dot{h})+  \ld \lb   G(t)+ \gam ^2 | \udp
|^2 + 2 \gam  \frac{d}{dt} \lb
   Q(t) + \eps(t)  \bar{\Theta}(u) \rb
   -2 \gamma  \dot{\eps}(t)  \bar{\Theta} (u)   \rb
 \\
 \hskip 2.cm  \le  |\udp|^2-\eps(t) \la   \nabla \Theta (u),  u-p \ra, \ea \]
or equivalently
\[
 w(t)+ \ld  G(t) + (\ld  \gamma ^2-1)  | \udp |^2  + \dot{\Gamma}_1 (t) \le
 -\eps(t) \la   \nabla \Theta (u),  u-p \ra  + 2 \ld \gamma  \dot{\eps}(t)  \bar{\Theta}
 (u),  \]
 where the function $\Gamma_1(\cdot)$ is defined by
 \[ \ba \Gamma_1(t)= \dot{h}(t) + \gam h(t) + 2 \ld \gam \lb Q(t) + \eps(t)  \bar{\Theta}(u(t))
 \rb \\
 \hskip 1.cm = \dot{h}(t) + \gam h(t) + \ld  \gamma  \lb  |\udp(t)|^2 +
 2   \lb  \phi(u(t)) + \la u(t)-p, Ap\ra \rb  + 2  \eps(t)  \bar{\Theta}(u(t)) \rb,   \ea \]
which by $Ap=-\nabla \phi(p)$ leads to (\ref{int1})\edem

\smallskip

We also need  the following variant of Gronwall's inequality,
see \cite{copeso} (Lemma 1).
\begin{lemma}\label{sup}
Let $\psi :[0,+\li) \to \R$ be absolutely continuous with
\[ \dot{\psi}(t) + \eps(t) \psi (t) \le \eps(t) g(t), \es \mbox{a.e.},\]
 where $g(t)$ is bounded and $\eps(t) \ge 0$ with $\eps \in
L^1_{\mbox{loc}}(\R_+)$. Then the function $\psi(t)$ is bounded and
if $\int_0^{+\li} \eps (\tau) d \tau = + \li$ we have $ \limsup_{ t \to
+\li} \psi (t) \le  \limsup_{ t \to +\li} g (t)$.
\end{lemma}
\ddem Let $\kappa_s:=\sup \{ g(t):  t \ge s\}$ so that
$\dot{\psi}(t)+ \eps(t) [ \psi(t)-\kappa_s] \le  0$ for $t \ge s$.
Multiplying by $\exp \lb \int_0^t \eps(\tau) d\tau \rb$ and
integrating over $[s,t]$ we get
$ \label{com1} [ \psi(t)-\kappa_s] \le [ \psi(s)-\kappa_s] \exp \lb
-\int_s^t \eps(\tau) d\tau \rb $.
It follows that $\psi(t)$ is bounded and, if $\int_0^{+\li}
\eps(\tau)d\tau =+\li$, by letting $t \to +\li$ in the above estimate
we obtain $\limsup_{t \to +\li} \psi (t) \le \kappa_s$. Letting $s \to
+\li$ in this last inequality yields $\limsup_{t \to +\li} \psi (t) \le \limsup_{t \to +\li}
g(t)$\edem

\smallskip

\ind  At once  we prove the boundedness of the trajectory given  by
(\ref{pb1}).
\begin{lemma} \label{yes}
Let   conditions   {\rm (H1)-(H4)} be satisfied with parameters such
that \bc $ \ld  \gamma^2
> 1   $, $\dot{\eps}(t) \le 0 $  for all $t\geq 0$, \ and \ $\eps (t) \to 0^+$ as $t \to
+\li$. \ec
  Then for any  solution $u$ of (\ref{pb1})  it holds that \\
 \ind   $u \in L^{\infty}([0,+\li),\Hs)$ and $\dot{u} \in L^{\infty}([0,+\li),\Hs)$. \\
 If in addition  $\int_0^{+\li} \eps(s) ds =
   +\li$, then the following properties are equivalent: \\
   \ind (i1) any weak cluster point of $u(t)$ for $t \to + \li$
   belongs to $S:=(\nabla \phi + A)^{-1}(0)$, \\
   \ind (i2)  $ \liminf_{ t \to +\li} \la \nabla \Theta(u_*),u(t)-u_*\ra \ge 0 $, \\
   \ind (i3) $u(t) \to u_*$ strongly,\\
   where $u_*$ is the
    unique solution of (\ref{pbv}).
  \end{lemma}
   \ddem Let us take $p = u_*$ in Lemma \ref{res2}, use the same notations  and   set \\
\ind $z(.):=\la \nabla \Theta(u), u-u_* \ra$,
 $ K(.)=\phi(u) - \la u-u_*, \nabla \phi(u_*)\ra $, \\
 so that  $  \Gamma_1=\dot{h} + \gam h + \ld  \gamma  \lb  |\udp |^2 +
 2   K   + 2  \eps(t)  \bar{\Theta}(u)
 \rb$.
 With these notations, by using
   $\dot{\eps} \le 0$ and $\bar{\Theta}(u) \ge 0$,  Lemma \ref{res2} gives
 \be \label{uni}  \ba
 \dis \ld  G+ w  + (\ld \gam ^2-1) | \udp  |^2
  + \dot{\Gamma}_1    \le - \eps  z . \ea \ee
For any real  values  $\nu >0$ and $\rho \in (0,1)$, by  definition
of $\Gamma_1$ we also have \[
  \nu  \Gamma_1-  z =  \lb -z+ \rho  \eta h + \nu 2 \ld \gam \eps
\bar{\Theta} (u) \rb  +   \lb \nu \gam -  \rho  \eta
 \rb h+ \nu  \ld  \gam |\udp |^2 + \nu  \dot{h}
 + 2   (\ld \gam \nu ) K ,
 \]
 which by (\ref{uni}) amounts to
 \be  \label{jit}  \ba
 \ld  G+ w  + (\ld \gam ^2-1) | \udp |^2
  + \dot{\Gamma}_1  + \nu \eps \Gamma_1 \\
   \le  \eps \lb -z+ \rho  \eta h + \nu 2 \ld \gam \eps
\bar{\Theta} (u) \rb + \eps  \lb \nu \gam -  \rho \eta
 \rb h   + \nu \eps \ld  \gam |\udp |^2 + \nu \eps \dot{h}
 + 2 \eps  (\ld \gam \nu ) K .
 \ea \ee
 By convexity of $\phi$ and  definition of $K$
 we  have
 $  K(t) \le \phi(u_*) + w(t)$,
 while by Young's inequality
  $ \dot{h} \le  (1/2)|\dot{u}|^2 + h$.
 As a result,  combining these last two inequalities  with   (\ref{jit}), we deduce
 that
\be \label{jit2}  \ba
 \dis  \ld  G+ ( 1- 2 \ld \gam \nu\eps  )w
  + [ \ld \gam ^2-1-\nu \eps (\ld \gamma+ 1/2)] | \udp  |^2
  + \dot{\Gamma}_1  + \nu \eps \Gamma_1  \\
 \hskip 1.cm \le   \eps \lb -z + \rho  \eta h  +  2
\ld  \gam \nu\eps  \bar{\Theta} (u) \rb   + \eps \lb \nu \gam + \nu
-  \rho  \eta
 \rb h +  2 \ld \gam \nu  \eps \phi(u_*).
 \ea \ee
Hence, for $\nu =(\rho  \eta)/(\gam+1)$, we equivalently obtain
 \be
 \label{rel1}  \ba
 \dis \ld  G +(1-2 \ld \gam \nu\eps )w
   + (\ld  \gam ^2-1- \nu \eps (\ld
  \gam+1/2)) | \udp  |^2
  + \dot{\Gamma}_1   + \nu \eps   \Gamma_1 \\
\hskip 2.cm \le   \eps \lb -z + \rho \eta h + 2\nu \ld \gam \eps
\bar{\Theta} (u) \rb  + 2  \ld \gam \nu \eps \phi(u_*).
 \ea \ee
  Moreover, since   $\ld  \gamma^2 -1 >0$ and  $\eps(t) \to 0$ as $t \to +\li$,
 for $t$ large enough, say   $t \ge
t_0$, we have  $ \ld  \gam ^2-1- \nu \eps(t) (\ld  \gam +1/2)  \ge
0$, $1-2 \ld \gam \nu\eps(t) \ge 0$, which by (\ref{rel1}) and
positivity of $G(\cdot)$ (and omitting the variable $t$) leads to
 \be \label{puis1}
  \ba
 \dis  \dot{\Gamma}_1   + \nu \eps   \Gamma_1
  \le   \eps \lb -z + \rho  \eta h  +  2 \ld  \gam \nu\eps
\bar{\Theta} (u) \rb  +  2  \ld \gam \nu \eps \phi(u_*).
 \ea \ee
 By convexity of $\Theta$, we know that
 \be \label{der1} z \ge \bar{\Theta}(u) -  \bar{\Theta}(u_*), \ee
 while by $\eta$-strong-monotonicity of $\nabla \Theta$ we have
 \be  \label{der2} z \ge \la \nabla \Theta(u_*),u-u_*\ra + 2 \eta h.  \ee
 From Young's inequality this latter inequality yields
\be  \label{der3} z \ge -(1/2 \eta) |  \nabla \Theta(u_*)|^2 + \eta
h  . \ee
Recalling that $\rho \in (0,1)$, and using (\ref{der1})
 and   (\ref{der3}), we obtain
 \be  \label{der4}  z \ge (1- \rho) ( \bar{\Theta}(u) -
 \bar{\Theta}(u_*))  -\rho(1/2 \eta) |  \nabla \Theta(u_*)|^2 +
\rho \eta h , \ee
 which by (\ref{puis1}) amounts to
  \be \label{puis2}
  \ba
 \dis  \dot{\Gamma}_1   + \nu \eps   \Gamma_1
  \le   \eps \lb (1- \rho)  \bar{\Theta}(u_*)
  + \rho(1/2 \eta) |  \nabla \Theta(u_*)|^2 \rb \\
\hskip 3.cm   - \eps \lb    (1- \rho)  -  2 \ld \gam \nu\eps
 \rb \bar{\Theta} (u) +  2 \ld \gam \nu  \eps \phi(u_*).
 \ea \ee
 Observing that $ (1- \rho)  -  2 \ld \gam \nu\eps(t) \ge 0$
  for $t \ge t_1$  large enough, we
 obtain
 \be \label{rel2}  \ba
 \dis  \dot{\Gamma}_1   + \nu \eps   \Gamma_1
  \le  \eps \lb (1- \rho)  \bar{\Theta}(u_*)
  + \rho(1/2 \eta) |  \nabla \Theta(u_*)|^2 \rb
  +  2 \ld \gam \nu  \eps  \phi(u_*).
 \ea \ee
 Applying Lemma \ref{sup} we deduce that $\Gamma_1$ is bounded.
 By  definition of $\Gamma_1$ and  convexity of $\phi$  (note that $\phi(u_*) \le K(t)$) we have
 \be \label{jl6}  \dot{h}(t) + \gam h(t) + \ld  \gamma  \lb  |\udp (t)|^2
   + 2 \phi(u_*)   \rb \le \Gamma_1(t). \ee
 Hence   there exists  a   positive constant  $C$ such that $t$ large enough, say $t \geq
 t_2$, yields
 \be \label{jil2}   \dot{h}(t) + \gam h(t) + \ld  \gamma   |\udp(t)|^2
 \le C, \ee
  which obviously  implies that  $h$ is bounded on $[0,+\li[$, and so is $u$:
  \be \label{jil3}
\sup_{ t \in [0,+ \li)} |u(t)| < +\li. \ee
 From $\dot{h}=\la \dot{u}, u-u_* \ra$, (\ref{jil2}) and  (\ref{jil3}) we
 also derive the boundedness of $\dot{u}$:
 \be \label{jil4}
\sup_{ t \in [0,+ \li)} |\dot{u}(t)| < +\li. \ee
Let us now prove the
equivalences. Note that  (i1) $\Rightarrow$ (i2) follows from
(\ref{pbv}), while (i3) $\Rightarrow$ (i1) is obvious. It remains to
prove  (i2) $\Rightarrow$ (i3). Let us assume that (i2) holds.
Again from (\ref{puis1}) and (\ref{der2}) we immediately have
 \be \dot{\Gamma}_1   + \nu \eps \Gamma_1
  \le  \eps \lb -  \la \nabla \Theta(u_*), u-u_*\ra    +  2 \ld  \gam \nu \eps
\bar{\Theta} (u) \rb  +  2 \ld \gam \nu  \eps  \phi(u_*)\ee
 which in light of  Lemma
\ref{sup} entails
 \[\ba  \limsup_{ t \to +\li} \Gamma_1 (t) \\
 \hskip 1cm  \le (1/\nu) \limsup_{ t \to +\li}
  [  -  \la \nabla \Theta(u_*),u(t)-u_*\ra    +  2 \ld  \gam \nu \eps(t)
\bar{\Theta} (u(t))  +  2 \ld \gam \nu  \phi(u_*) ] \\
  \hskip 1cm =(1/\nu) [   2 (\ld \gam \nu) \phi(u_*)
   - \liminf_{ t \to +\li} \la \nabla \Theta(u_*),u(t)-u_*\ra    ] .
\ea \]
It is then immediate  that   (i2) leads to
 $  \limsup_{ t \to +\li} \Gamma_1 (t) \le 2 \ld \gam \phi(u_*)$,
 which in light of (\ref{jl6}) entails
 $ \limsup_{ t \to +\li}  \lb \dot{h}(t)  + \gam    h(t) \rb \le 0$.
 Then, it is easily checked (one can apply Lemma \ref{sup}) that this inequality implies
   $\limsup_{ t \to +\li} h(t)=0$, namely (i3)\edem

\smallskip
\ind We  are now in position to prove  convergence in norm of the
trajectory of the solution $u$ of (\ref{pb1}) towards a special zero
of $\nabla \phi + A$.

{\bf Proof of Theorem \ref{imp}}. From Lemma \ref{yes} we know that
 \be \label{est1} \mbox{   $u \in L^{\infty}([0,+\li);\Hs)$ and
$\dot{u} \in L^{\infty}([0,+\li);\Hs)$}. \ee
 By choosing $p= u_*$ in  Lemma \ref{res2},  keeping the same
 notations, and setting $ K(t) = \phi(u(t)) - \la u(t)-u_*, \nabla \phi(u_*)\ra
 $, we recall that
\begin{center}
$  \Gamma_1(t)=\dot{h}(t) + \gam h(t) + \ld  \gamma  \lb  |\udp(t)|^2 +
 2   K(t)   + 2  \eps(t)  \bar{\Theta}(u(t)) \rb$,
\end{center}
 and we also have
  \be \label{ok} \ld G(t) + w(t) + (\ld  \gam ^2-1)
   | \udp(t) |^2
+  \dot{\Gamma}_1 (t)  \le - \eps(t) \la \nabla \Theta (u(t)), u(t)-u_* \ra .
\ee
The remainder of the proof will be divided into two cases:

\ind Case 1.   Assume that there exists a time $t_0$ such that the
function $\Gamma_1(\cdot)$ is nonincreasing on $[t_0,+ \li)$. Then recalling
that $\Gamma_1$ is lower bounded (thanks to its definition and using
(\ref{est1})),  we
  deduce that $ \Gamma_1 (t)$ converges as $t \to +\li$ to   some real value
  $\al$. This leads to
   $\lim_{t \to +\li} \dot{\Gamma}_1  (t)= 0 $. Indeed, for any reals $t$ and $s$ verifying
   $t
  \ge s \ge t_0$, we have
  $ \int_{s}^t |\dot{\Gamma}_1 (\tau)|d\tau= -\int_{s}^t \dot{\Gamma}_1 (\tau)d\tau=
  \Gamma_1(t) -\Gamma_1(s)$,
  hence, passing to the limit as $t \to +\li$, we get
   $ \int_{s}^{+\li} |\dot{\Gamma}_1 (\tau)|d\tau= \al - \Gamma_1(s) < +\infty$, and,
  by Lipschitz continuity of $|\dot{\Gamma}_1 (\cdot)|$, we conclude that $\dot{\Gamma}_1 (t) \to 0$ as $t \to +\li$.
  Therefore, by  (\ref{ok}),  the
boundedness of $u(\cdot)$, together with $\lambda \gam^2 -1>0$ and $\eps(t)
\to 0$ as $t \to +\li$, we deduce that
$ \lim_{t \to  +\li} w(t) =\lim_{t \to  +\li} G(t)=\lim_{t \to
+\li}| \dot{u}(t)| = 0 .$
It follows that  $ \lim_{t \to  +\li} w(t) =0$ and $\lim_{t \to +\li}
|Au(t) - Au_*|=0$. Hence,  in light of Lemma \ref{wcp}, we observe
that  any weak cluster point of $u$ belongs to $S$,
  which   by  Lemma \ref{yes} shows  that $u(t) \to u_*$ strongly.

\ind Case 2.  Assume that there exists an increasing  sequence
$(l_k) \subset [0,+\li)$ such that
 $\lim_{k \to \li}
l_k=+\li$ and $\dot{\Gamma}_1 (l_k) >0$ for all $k \ge 0$. This allows us to
introduce  the  mapping $s: [0,+\li) \to [0,+\li)$ defined  for $t$
large enough by
   \[ s(t) = \sup \left \{ \beta \le t \es | \es \dot{\Gamma}_1 (\beta) >0 \right \} ,\]
 and we establish  that $\limsup_{t \to +\li}  \Gamma_1 (t) \le
   0$.  Let $(t_n) \subset [0,+\li)$ be a sequence such that $\lim_{n \to +\li}
t_n=+\li$, and set $s_n = s(t_n)$.  It is obviously seen that
$(s_n)$ is a nondecreasing sequence such that $s_n \to +\li$ as
$n \to +\li$, $(\dot{\Gamma}_1 (s_n)) \subset  [0,+\li)$, while  by
(\ref{ok}) we have
  \be \label{ok2} \ba
 \ld G(s_n)+ w(s_n) + (\lambda \gam ^2-1)
   | \udp(s_n) |^2  \le  - \eps (s_n) \la \nabla \Theta(u(s_n)), u(s_n)-u_* \ra . \ea  \ee
 From this last inequality,  keeping in mind that  $\ld  \gam^2 -1>0$ and $\eps(s_n) \to 0$ as $n \to
\li$, we obviously deduce that
 \be \label{vu1}
 \lim_{n \to +\li} | \udp(s_n) |= \lim_{n \to +\li} w(s_n)=\lim_{n \to +\li} | Au(s_n)-Au_*|=0. \ee
By using $\eps(s_n) \in  (0,+\li)$, we
additionally obtain  \be \label{vu2} \la \nabla \Theta (u(s_n)),
u(s_n)-u_* \ra \le 0. \ee
 Clearly, (\ref{vu1}) entails that any weak-cluster point of
 $u(s_n)$ belongs to $S$ (according to Lemma \ref{wcp}),
  while by (\ref{vu2}) and using the
 $\eta$-strong monotonicity
 of $\nabla \Theta$, we obtain
 \[  \eta |u(s_n)-u_*|^2 \le \la \nabla {\Theta} (u(s_n)) - \nabla \Theta (u_*),  u(s_n) -
 u_* \ra
 \le - \la \nabla \Theta (u_*),  u(s_n) -
 u_* \ra. \]
Hence
\[  \limsup_{n \to +\li} \lb \eta |u(s_n)-u_*|^2 \rb
 \le - \liminf_{n \to +\li} \la \nabla \Theta (u_*),  u(s_n) -
 u_* \ra \le 0, \]

\noindent this last inequality  being obtained by passing to the limit on a weak convergent subsequence of $u(s_n)$, using the fact that
 any weak-cluster point of
 $u(s_n)$ belongs to $S$ and the definition (\ref{pbv}) of $u_*$.
 Hence $u(s_n) \to u_*$ strongly, which together with (\ref{vu1}) and the definition of
  $\Gamma_1$ implies that
  \be \label{vu3} \lim_{n \to +\li}  \Gamma_1 (s_n)=2 \ld \gam \phi(u_*). \ee
   In addition, using the definition of $s(t_n)$
 we clearly have $s_n \le t_n$ and $\dot{\Gamma}_1(t) \le 0$ for $t \in
 ]s(t_n), t_n[ $ when $s(t_n) \neq t_n$.
 Let us write
 $ \Gamma_1(t_n) =  \Gamma_1 (s_n)+ \int_{s_n}^{t_n } \dot{\Gamma}_1(t) dt
 \le \Gamma_1 (s_n)$,
   which by (\ref{vu3})  yields
 $  \limsup_{n \to \li}  \Gamma_1 (t_n) \le
   2 \ld \gam \phi(u_*)$.
   By convexity of $\phi$
   $ \dot{h}(t) + \gam h(t) + 2 \ld \gam \phi(u_*) \le \Gamma_1(t)$,
   which yields
   $  \limsup_{n \to \li} \lb  \dot{h}(t_n) + \gam h(t_n) \rb \le
   0$.
   This being true for any sequence $(t_n)$ that converges to infinity, we
   conclude that
   $ \limsup_{t \to +\li} (\dot{h}(t) + \gam h(t)) \le 0$,
   which immediately leads to $u(t) \to u_*$ strongly \edem

\medskip

The case $\Theta (v)  = \frac{1}{2} |v|^2$,  which corresponds to the classical Tikhonov regularization method, has been the object of active research, especially in the case of  first order  dynamical systems
governed by maximal monotone operators,
see \cite{copeso} for some advanced results and references. Concerning second order dynamics,  the potential case has been considered in \cite{acz}.
Let us now state an extension of this last result to the non-potential case.

\begin{corollary} \label{corvis}    Let us suppose that   {\rm (H1)-(H3) } hold with  $S := ( \nabla
\phi + A)^{-1}(0) \neq \emptyset$ and with $  \ld \gamma^2  > 1$.
Let us assume moreover that $\eps (t)$ is decreasing and converges slowly to zero as $t
\to +\li$ in the following sense:

    \be
   \int_0^{+\li} \eps(s) ds = +\li.
    \ee

\noindent Then, for any  $u_0$ and $v_0$ in $\Hs$,
there exists a unique  solution $u \in C^2([0,+\li); \Hs)$ of

\be \label{Tikh}
  \uds (t) +\gam \udp(t)  + \nabla \phi(u(t))  +   A(u(t)) + \eps(t) u(t) = 0, \ee

\noindent with $u(0)=u_0$ and $\dot{u}(0)=v_0$, which satisfies
$u(t) \to u_*$ strongly in  $\Hs$   as $t \to +\li$, where $u_*$ is
the element of minimal norm of the closed convex set $S$.

\end{corollary}

 \section{\large Examples}
The following examples  aim at illustrating our abstract results. Each of them requires further studies which are out of the scope of the present article.

\subsection{Constrained optimization}
 Let $C$ be a closed convex subset of $\Hs$, and $g: \Hs
\to \R$  a convex continuously differentiable function, whose gradient $ \nabla g$ is
Lipschitz continuous over $\Hs$ with Lipschitz constant $L$. Let us consider the  constrained optimization problem
\be \label{gp1}
\min \left\{ g(v) :  v\in C            \right\}.
\ee
Because of their direct linking with numerical gradient methods, we are interested in the study of
smooth  dissipative dynamical systems
whose trajectories asymptotically converge towards minimizers of (\ref{gp1}).
Considering second-order dynamics aims at making the algorithms faster.

\noindent Let us denote by $S$ the set of solutions of (\ref{gp1})
and suppose that $S \neq \emptyset$.
 A first-order optimality condition for (\ref{gp1}) is given by
 the following inclusion

\be \label{gp2} \nabla g(u) + N_C (u) \ni 0, \ee

\noindent where $N_C(u)$ is the outward normal cone to $C$ at $u\in
C$. Observe that (\ref{gp2}) is equivalent to
\be \label{gp3}
u - P_C[u-\mu \nabla g(u) ]= 0,
\ee

\noindent where $P_C$ is the projection operator on $C$ and $\mu$ is a positive parameter.
 Whence we considerer the following gradient-projection dynamic
\be \label{gp4}
\uds (t) + \gam \udp(t)  +  u(t)- P_C[u(t)-\mu \nabla g(u(t)) ]= 0.
\ee

\noindent This  continuous  dynamical system has been first introduced by Antipin \cite {ant} and
 further studied by Attouch-Alvarez \cite {alat2}.   An extended study of the corresponding first-order system
has been achieved in \cite{Bol1}, \cite{Bol2}. Let us show that the
operator \be \label{gp5} Av: = v - P_C(v-\mu \nabla g(v)) \ee
\noindent  is $(1/2)$-cocoercive for $0 < \mu < (2/L)$. To that end,
let us notice that $A = I - T$ where $Tv = P_C(v-\mu \nabla g(v))$.
According to Corollary \ref{cor2}, and noticing that $P_C$ is a
contraction, we just need to verify that the mapping $v \mapsto
v-\mu \nabla g(v)$ is a contraction too. Given arbitrary $u\in \Hs$,
$v\in \Hs$, we have

\[ \ba  | (u-\mu \nabla g(u)) - (v-\mu \nabla g(v))|^2 \\
\hskip 1.cm =   |u - v |^2  -2\mu \left\langle \nabla g(u) - \nabla
g(v), u -v \right\rangle + {\mu}^2 |\nabla g(u) - \nabla g(v)|^2.
 \ea
 \]

\noindent Thus proving that $v \mapsto v-\mu \nabla g(v)$ is a contraction, is equivalent
 to prove that
%
$\left\langle \nabla g(u) - \nabla g(v), u -v \right\rangle  \geq
\frac{\mu}{2}|\nabla g(u) - \nabla g(v)|^2$,
%
\noindent which amounts to say that $\nabla g $ is
$(\mu/2)$-cocoercive. Clearly, this property forces $\nabla g $ to
be $(2 / \mu )$-Lipschitz continuous. Indeed, as a striking property,
the  converse statement holds true for the gradient of a convex
function, that's the Baillon-Haddad theorem, see \cite{BaHa}. Hence
$\nabla g $ being Lipschitz continous with Lipschitz constant $L
\leq (2 / \mu)$ forces $\nabla g $ to be $(\mu/2) $-cocoercive.

Hence  (\ref{gp4}) falls in our
 setting.   Applying Theorem \ref {wac}
with $\lambda = 1/2$, we deduce that, for $\gamma > \sqrt{2}$ and $0
< \mu \leq  (2/L)$,  \  every trajectory of (\ref{gp4}) weakly
converges to an element of $S$, which is a minimizer of problem
(\ref{gp1}).

Moreover, when considering the Tikhonov-like  regularized dynamic with $\eps(t) \to 0$ as $t$ goes to $+\infty$,

\be \label{gp6}
  \uds (t) +\gam \udp(t)  +  u(t)- P_C[u(t)-\mu \nabla g (u(t)) ] +   \eps(t) \nabla \Theta (u(t)) = 0,
\ee

\noindent by a direct application of Theorem \ref{imp}, under the
assumption  \  $ \int_0^{+\li} \eps(s) ds = +\li$, \ one obtains
that the  trajectories  of (\ref {gp6}) strongly converge to the
unique
 solution $\bar{u}\in S$ verifying
 \be
 \Theta (\bar{u})= \inf_{ \mbox{$v\in S$}} \Theta (v).
 \ee

\noindent  \textbf{Remark} Inertial dynamical systems, like the "heavy ball with friction dynamical system" have been first introduced by B. Polyak \cite{Po} in the realm of optimization. Since then, an
abundant literature has been devoted to this subject. Rich connections between second order dissipative dynamical systems in their respective continuous and discretized forms  have been put to the fore, and so provide new
algorithms  together with a deeper insight into their  convergence analysis,
 see \cite{alv}, \cite{alat},  \cite{alat2}, \cite{alat3} , \cite{alat4}, \cite{acb},
 \cite{acz}, \cite{cab}, \cite{flam}, \cite{main1}.

\subsection{Coupled systems and dynamical games }

Throughout this section we make the following standing assumptions:

\begin{itemize}
 \item $\mathcal H = {\mathcal X}_1 \times {\mathcal X}_2$ is  the cartesian product of two Hilbert spaces
 equipped with norms $|.|_{{\mathcal X}_1}$ and $|.|_{{\mathcal X}_2}$,  while
$x=(x_1,x_2)$, with $x_1 \in {\mathcal X}_1$ and $x_2 \in {\mathcal
X}_2$, stands for any element in $\Hs$;
\item  $\phi(x) = f_1(x_1) + f_2(x_2) + \Phi (x_1,x_2) $,
where   $f_1: {\mathcal X}_1 \to  \R$,  $f_2: {\mathcal X}_2  \to  \R$
are smooth convex functions, $\Phi : {\mathcal X}_1 \times {\mathcal X}_2 \to \R$
is a smooth convex coupling function;
\item  $A = ({\nabla}_{x_1}\mathcal L, -{\nabla}_{x_2}\mathcal L)$ is the maximal monotone operator which
 is attached to a smooth
convex-concave function $\mathcal L:  {\mathcal X}_1 \times {\mathcal X}_2 \to \R$. The operator
$A$ is assumed to
be $\lambda$-cocoercive for some $\lambda >0$.
\end{itemize}

\noindent A typical example of coupling function $\Phi$ (see \cite{abrs}) is given by

\begin{center}
$\Phi (x) = \frac{1}{2}|L_1x_1 - L_2x_2|^2_{\mathcal Z}$,
\end{center}

\noindent where $L_1 \in L({\mathcal X}_1, \mathcal Z)$ and  $L_2
\in L({\mathcal X}_2, \mathcal Z)$ are linear continuous operators
acting respectively from ${\mathcal X}_1$  and ${\mathcal X}_2$ into
a third Hilbert space $\mathcal Z$ with norm $|. |_{\mathcal Z}$.

\noindent A smooth convex-concave function $\mathcal L$, such that
the associated maximal monotone
 operator $A$ is $\lambda$-cocoercive, can be obtained in a systematic way by using  epi-hypo
 Moreau-Yosida regularization, as described below. Given
 a general closed convex-convave function $\mathcal L:  {\mathcal X}_1 \times {\mathcal X}_2 \to \R$ (see \cite{Rock}) and  a positive parameter
  $\lambda$, we consider  $\mathcal L_{\lambda} :  {\mathcal X}_1 \times {\mathcal X}_2 \to \R$ defined for any $
  (x_1,x_2) \in {\mathcal X}_1 \times {\mathcal X}_2$ by

\begin{center}
$ \dis   {\mathcal L}_{\lambda} (x_1,x_2) = \min_ {u_1 \in {\mathcal
X}_1} \max_ {u_2 \in {\mathcal X}_2} \left\{  \mathcal L(u_1,u_2) +
\frac{1 }{2\lambda } |x_1 -u_1|_{{\mathcal X}_1}^2 -  \frac{1
}{2\lambda } |x_2 -u_2|_{{\mathcal X}_2}^2       \right\}$.
\end{center}

\noindent Then $ {\mathcal L}_{\lambda} $ is a smooth convex-concave
function whose associated maximal monotone operator
$({\nabla}_{x_1}{\mathcal L}_{\lambda}, -{\nabla}_{x_2}{\mathcal
L}_{\lambda})$ is precisely the Yosida approximation $A_{\lambda}$
of the operator $A = ({\partial}_{x_1}\mathcal L,
-{\partial}_{x_2}\mathcal L)$,
 i.e., $A_{\lambda} = ({\nabla}_{x_1} \mathcal L_{\lambda}, -{\nabla}_{x_2}\mathcal L_{\lambda})$ (see \cite{aaw} for further details).

\noindent In this setting,  with $ u(t) = (x_1(t),x_2(t)) $ system (\ref{pb00}) becomes

\begin{equation}\label{couplsys1}
\hskip -0.5cm \left\{
\begin{array}{l}
\ddot{x}_1(t) +  \gamma \dot{x}_1(t)
 +   \nabla f_1(x_1(t)) + {\nabla}_{x_1}\Phi (x_1 (t),x_2 (t))
   +  {\nabla}_{x_1}\mathcal L  (x_1 (t),x_2 (t)) =  0,\\
\rule{0pt}{25pt}
\ddot{x}_2(t) +  \gamma \dot{x}_2(t)   +   \nabla f_2(x_2(t)) + {\nabla}_{x_2}\Phi (x_1 (t),x_2 (t)) -  {\nabla}_{x_2}\mathcal L  (x_1 (t),x_2 (t)) =  0.
\end{array}\right.
\end{equation}

\noindent As a straight application of Theorem  \ref {wac}, assuming
relation $\lambda {\gamma}^2 > 1$ holds, one obtains $x(t)=
(x_1(t),x_2(t))\longrightarrow
x_{\infty}=(x_{1,\infty},x_{2,\infty})$ weakly in  $\mathcal H$ as
$t$ goes to $+\infty$, where $(x_{1,\infty},x_{2,\infty})$ is
solution of the coupled system

\begin{equation}\label{couplsys2}
\left\{
\begin{array}{l}
  \nabla f_1(x_1) + {\nabla}_{x_1}\Phi (x_1 ,x_2)  +  {\nabla}_{x_1}\mathcal L  (x_1,x_2) =  0,\\
\rule{0pt}{15pt}
 \nabla f_2(x_2) + {\nabla}_{x_2}\Phi (x_1,x_2) -  {\nabla}_{x_2}\mathcal L  (x_1,x_2) =  0.
\end{array}\right.
\end{equation}

\noindent Structured systems such as    (\ref{couplsys2})  contain
both potential and non-potential terms whose antagonistic effects
are often present in decision sciences and physics.

 In game theory, (\ref{couplsys2}) describes Nash equilibria of the normal form game with two players $1$,  $2$ whose
static loss functions are respectively given by
\begin{equation}\label{couplsys3}
\left\{
\begin{array}{l}
F_1:(x_1,x_2)\in {\mathcal X}_1 \times {\mathcal X}_2 \rightarrow F_1(x_1,x_2)
 = f_1(x_1) +  \Phi(x_1,x_2)  +  \mathcal L  (x_1,x_2),\\
\rule{0pt}{15pt} F_2:(x_1,x_2)\in {\mathcal X}_1 \times {\mathcal
X}_2 \rightarrow F_2(x_1,x_2) = f_2(x_2) +  \Phi(x_1,x_2) - \mathcal
L  (x_1,x_2),
\end{array}\right.
\end{equation}

\noindent where the  $f_i (.)$ represent the individual convex
payoffs of the players, $\Phi(.,.)$ is their joint convex payoff,
and $\mathcal L $ is a convex-concave payoff with  zero-sum rule.
The case $\mathcal L = 0$ corresponds to a potential team game (see
\cite{MS}), while case $\Phi = 0$ corresponds to a  non cooperative
zero-sum game. Game (\ref{couplsys3}) involves both cooperative and
non cooperative aspects.

\noindent A central question in game theory, decision sciences and economics is to describe realistic dynamics which converge to Nash equilibria, see \cite{HS} and the references herein.
Indeed, implicit time discretization of  dynamical system (\ref{couplsys1}) leads to the following best response dynamics (players $1$ and $2$ play alternatively) with inertia and ``costs to change''
\begin{center}
 $(x_{1,k+1}, x_{2,k+1})\rightarrow (x_{1,k+2}, x_{2,k+1})\rightarrow (x_{1,k+2}, x_{2,k+2})$;
\end{center}
\begin{equation}\label{couplsys4}
\left\{
  \begin{array}{l}
 x_{1,k+2} = \argmin_{\xi\in {\mathcal X}_1}
    \{  f_1(\xi)+ \Phi(\xi,x_{2,k+1})  +
    \mathcal L  (\xi,x_{2,k+1}) \\
  \hskip 4.cm   + \frac{1}{2{\alpha}_k} \parallel \xi -\left( x_{1,k+1}+
     {\beta}_k (x_{1,k+1} - x_{1,k}) \right)\parallel^{2}_{{\mathcal X}_1}\},  \\
\rule{0pt}{15pt}
 x_{2,k+2}= \argmin_ {\eta\in {\mathcal X}_2}
    \{f_2(\eta)+  \Phi( x_{1,k+2},\eta) - \mathcal L  (x_{1,k+2},\eta) \\
 \hskip 4.cm    + \frac{1}{2{\nu}_k} \parallel \eta - \left( x_{2,k+1}+  {\beta}_k (x_{2,k+1} - x_{2,k}) \right)\parallel^{2}_{{\mathcal X}_2}
     \}.\\
   \end{array}\right.
\end{equation}
%
 The  terms $\parallel \xi - x_{1,k+1}-  {\beta}_k (x_{1,k+1} -
x_{1,k})\parallel^{2}$ and $\parallel \eta - x_{2,k+1} - {\beta}_k
(x_{2,k+1} - x_{2,k})\parallel^{2}$ naturally come into play when
discretizing first and second order time derivatives, see
\cite{alv}, \cite{alat}. In decision sciences, they reflect some
aspects  of  agents behaviors like anchoring and inertia,  and can
be interpreted as ``low local costs to change'', see
  \cite{abrs}, \cite{as1}. The various
parameters ${\alpha}_k, {\nu}_k, {\beta}_k$ model  adaptive,
learning  abilities of the agents as well as their reactivity or
resistance to change.

\noindent Continuous dynamical system (\ref{couplsys1}) together with its convergence properties offers
a valuable guideline to  study  discrete dynamic (\ref{couplsys4}). Clearly, this  requires  further
study beyond the limit of the present paper.

\subsection{Further perspectives}
It is natural  to relax the regularity assumptions on the potential
operator $\nabla \phi$ and the maximal monotone operator $A$, and to
wonder  whether our asymptotic convergence results   remain valid in
respect of the more general differential inclusion
 \be \label{mulva}  \ddot{u}(t) + \gamma \dot{u}(t)+ A u(t) +  \der \phi(u(t)) \ni 0. \ee
 \noindent Let us give   some precisions (regarding cases when either $A$ or $\der\phi$ is single-valued): \\
 \ind  1.   The non-smooth
potential case (obtained by taking $A=0$) , that is
 \be \label{mulva1}  \ddot{u}(t) + \gamma \dot{u}(t) +  \der \phi(u(t)) \ni 0, \ee
\noindent has been  considered in the finite dimensional case in \cite{acb}. Trajectories of
(\ref{mulva1}) naturally exhibit elastic shocks,  which makes this
system play an important role in unilateral
 mechanics, see \cite{Sch}.
Taking $\phi$ equal to the indicator function of a closed convex set $C$ yields  (damped) billiard dynamics.
 Given a non-smooth closed convex
potential $\phi : \Hs \to \R \cup  \{+\li\}$  and a
$\lambda$-cocoercive operator $A$, a natural approach to
(\ref{mulva}) (in the same   lines as \cite{acb}, \cite{Sch})  would consist in
 approximating $\phi$ by
 smooth convex functions $\phi_n$ (this is always possible
  for example via Moreau-Yosida approximation), which  leads  to study
  \be \label{mulva2}  \ddot{u}_n(t) + \gamma \dot{u}_n(t) +  Au_n(t) +  \nabla {\phi}_n(u_n (t)) = 0. \ee

Note that our basic assumption $\lambda {\gamma}^2 > 1$, which does
not concern the potential part of the operator, is not affected by
this operation. Thus, (\ref{mulva2}) falls in our setting. Then, one
has to establish estimations on the sequence   $(u_n)$,  pass to the
limit on these estimations as $n$ goes to infinity, and derive
conclusions on the asymptotic behavior of trajectories of
(\ref{mulva}). As it has been done in \cite{alat3} and \cite{acb},
this program  would allow to consider the asymptotic behavior of
damped hyperbolic equations (with non-isolated equilibria) combining
potential with non potential effects.

2. The study of (\ref{mulva})  with  a general maximal monotone
operator $A$  leads to considerable difficulties. In that case,
another type of regularization or relaxation method can be used. It
relies on the following remark. Let  $\lambda >0$ and $\gamma
>0$ (being fixed) such that $\lambda {\gamma}^2 > 1$,
 $A_{\lambda}$  the Yosida approximation of  $A$, and  consider the
 equation
\be  \label{pers0}
 \uds (t) + \gam \udp(t)   + \nabla \phi(u(t))+   A_{\lambda}(u(t))  =
 0.
 \ee
Equation (\ref{pers0}) enters  our setting (as $A_{\lambda}$ is
$\lambda$-cocoercive) and it can be equivalently  rewritten as

\begin{equation}\label{pers}
\left\{
\begin{array}{l}
\uds (t) + \gam \udp(t)   + \nabla \phi(u(t))+  \frac{1}{\lambda}(u(t) - v(t)) =  0,\\
\rule{0pt}{15pt}
Av(t) +  \frac{1}{\lambda}(v(t) - u(t))  \ni 0.
\end{array}\right.
\end{equation}

This naturally  suggests that some of our results can be extended to coupled systems of the following type

 \begin{equation}\label{pers1}
\left\{
\begin{array}{l}
\uds (t) + \gam \udp(t)   + \nabla \phi(u(t))+  \frac{1}{\lambda}(u(t) - v(t)) =  0,\\
\rule{0pt}{15pt}
\dot{v}(t) + Av(t) +  \frac{1}{\lambda}(v(t) - u(t))  \ni 0.
\end{array}\right.
\end{equation}

\noindent Our  approach, which is based on Liapunov methods,  (\cite{Har}, \cite{HarCaz}, \cite{HalRau}), can be adapted to such nonlinear systems.
 In the  case of linear hyperbolic PDE's, it would be worthwhile   to compare it with  spectral analysis methods,
  see for example \cite{EZ},  \cite{Ga}.

 3. In view of further applications,  it would be interesting   to consider other types of dissipation phenomena:
For example, dry friction  most likely produces stabilization within
finite time, a desirable feature of human decision processes, see
\cite{aac}. Geometrical damping is of great importance  for
numerical optimization, and modeling of non elastic shocks in
mechanics, see \cite{alat4}.


\footnotesize
\begin{thebibliography}{45}


\bibitem{aac} S. Adly, H. Attouch, A. Cabot, Finite time stabilization of nonlinear oscillators subject to dry friction, Nonsmooth mechanics and analysis,
 {\it Adv. Mech. Math}  \textbf{12} (2006), pp. 289--304.

\bibitem{alv} F. Alvarez, On the minimizing property of a second
order dissipative system in Hilbert space, {\it SIAM J. Control
Optim.}  \textbf{38} (4) (2000), pp. 1102--1119.

\bibitem{alv2} F. Alvarez, Weak convergence of a relaxed and inertial hybrid projection-proximal
point algorithm for maximal monotone operators in Hilbert space,
 {\it SIAM J. Optim.}  \textbf{14} (3) (2004), pp. 773--782.
%
\bibitem{alat} F. Alvarez, H. Attouch, An inertial proximal method for monotone operators via discretization of a nonlinear oscillator
with damping, {\it Set Valued Analysis}  \textbf{9} (2001), pp. 3--11.
%
\bibitem{alat2} F. Alvarez, H. Attouch,  The heavy ball with friction dynamical system for
convex constrained minimization problems. Optimization (Namur, 1998), pp. 25--35,
{\it Lecture Notes in Econom. and Math. Systems}  \textbf{481}, Springer, Berlin, (2000).
%
\bibitem{alat3} F. Alvarez, H. Attouch, Convergence and asymptotic stabilization for some damped
hyperbolic equations with non-isolated equilibria, {\it ESAIM Control Optim. Calc. Var.}
\textbf{6}, (2001), pp. 539--552.

\bibitem{alat4} F. Alvarez, H. Attouch, J. Bolte, P. Redont, A second-order gradient-like dissipative dynamical system with Hessian-driven damping.
 Application to optimization and mechanics,
    {\it  J. Math. Pures Appl.}  \textbf{81}, (2002), pp.  747--779.

\bibitem{ant} A.S. Antipin, Minimization of convex functions on convex sets by means of differential equations (in Russian), Differ. Uravn. 30 (9) (1994) 1475--1486;
English translation: {\it Differential Equations} \textbf{30} (9) (1994) pp. 1365--1375.

\bibitem{aaw} H. Attouch, D. Aze, R. Wets, Convergence of convex-concave saddle functions: Applications to convex programming
and mechanics, {\it Ann. Inst.  Poincar\'e}  \textbf{5} (6), (1988), pp. 537--572.

\bibitem{abrs} H. Attouch, J. Bolte, P. Redont, A. Soubeyran, Alternating proximal algorithms
for weakly coupled minimization problems. Applications to dynamical games and PDE's,
  {\it J. of Convex Analysis}  \textbf{15} (3),  (2008),  pp.  485--506.

\bibitem{acb} H. Attouch, A.  Cabot, P. Redont,  The dynamics of elastic shocks via epigraphical
 regularization of a differential inclusion. Barrier and penalty approximations,
  {\it Adv. Math. Sci. Appl.}  \textbf{12} (1),  (2002),  pp.  273--306.

\bibitem{atcom} H. Attouch, R. Cominetti, A dynamical approach to
convex minimization coupling approximation with the steepest descent
method, {\it J. Diff. Equat.} \textbf{128} (1996), pp. 519--540.

\bibitem{acz} H. Attouch, M.-O. Czarnecki,  Asymptotic control and stabilization of nonlinear
 oscillators with non-isolated equilibria,  {\it J. Diff. Equat.}
  \textbf{179 }(1), (2002), pp. 278--310.

\bibitem{agr} H. Attouch, X.  Goudou, P. Redont,  The heavy ball with friction method: The continuous dynamical system. Global
exploration of local minima by asymptotic analysis of a dissipative dynamical system.
regularization of a differential inclusion. Barrier and penalty approximations.
{\it Commun. Contemp. Math.}  \textbf{1},  (2000), pp. 1--34.

\bibitem{as1} H. Attouch, A. Soubeyran, Inertia and reactivity in decision making as cognitive variational inequalities,
  {\it J. Convex. Anal.}
  \textbf{13} (2), (2006), pp. 207--224.

\bibitem{BaHa} J.-B. Baillon and G. Haddad, Quelques propri\'etes des op\'erateurs angles-born\'es et n-cycliquement monotones,
      {\it  Israel J. Math.} \textbf{26}, (1977), pp. 137--150.

\bibitem{Bol1} J. Bolte, Continuous gradient projection method in Hilbert spaces,
      {\it  J. Optim. Theory Appl.} \textbf{119} (2), (2003), pp. 235--259.

\bibitem{Bol2} J. Bolte, M. Teboulle, Barrier operators and associated gradient-like dynamical systems
for constrained minimization problems,
{\it   SIAM J. Control Optim.} \textbf{42} (4), (2003), pp. 1266--1292.

\bibitem{Bre} H. Br\'ezis, Op\'erateurs maximaux monotones et semi-groupes de contractions dans les espaces de
Hilbert, North-Holland, Mathematical Studies, 1973.
%
\bibitem{cab}  A. Cabot, Inertial gradient-like dynamical system
controlled by a stabilizing term, {\it J. Optim. Theory Appl.} \textbf{120} (2), (2004),
pp. 275--303.

\bibitem{HarCaz}  T. Cazenave, A. Haraux,  An introduction to semilinear evolution equations,
 Oxford Lecture Series in Mathematics and its Applications  \textbf{13}, 1998.
%
\bibitem{combh} P.L. Combettes, S.A. Hirstoaga, Visco-penalization
of the sum of two operators, {\it Nonlinear Analysis: TMA}  \textbf{69} (2),
(2008),  pp. 579--591.
%
\bibitem{copeso} R. Cominetti, J. Peypouquet, S. Sorin, Strong
asymptotic convergence of evolution equations governed by maximal
monotone operators with Tikhonov regularization, {\it J. Diff.
Equat.} \textbf{ 245} (2008), pp. 3753--3763.

\bibitem{EZ} S. Ervedoza, E. Zuazua, Uniformly exponentially stable approximations for a class of damped systems,
 {\it J. Math. Pures et Appl.}  \textbf{91},
(2009),  pp. 20--48.

\bibitem{flam} S.D. Flam, J. Morgan,
Newtonian mechanics and Nash play, {\it International Game Theory Review}  \textbf{6}  (2), (2004), pp. 181--194.

\bibitem{Ga}  I. Gallagher,    Asymptotics of the solutions of hyperbolic equations with a skew-symmetric perturbation,
{\it Journal of Differential Equations} \textbf{150},  (1998) , pp.~ 363--384.

\bibitem{HalRau} J. K. Hale,    G. Raugel,
Convergence in gradient-like systems with applications to PDE,
Zeitschrift für Angewandte Mathematik und Physik, Birkhauser Verlag  Basel,
 \textbf{43}  (1), (1992), pp.~ 63--125.

 \bibitem{Har}  A. Haraux,  {\it Syst\`emes dynamiques dissipatifs et applications},  Masson, RMA \textbf{17}, 1991.

\bibitem{HS} J. Hofbauer and S. Sorin,  Best response dynamics for continuous zero-sum games,
{\it Discrete and Continuous Dynamical Systems, series B} \textbf{6 } (1), (2006), pp.~215--224.

\bibitem{main1} P.E. Maing\'e, Regularized and inertial
algorithms for common fixed points of nonlinear operators, {\it J.
Math. Anal. Appl.}  \textbf{344} (2008), pp. 876--887.

\bibitem{MS}  D. Monderer, L. S. Shapley, Potential Games, {\it Games
Econ. Behav.} \textbf{14} (1996), pp.~124--143.

\bibitem{Op} Z. Opial, Weak convergence of the sequence of successive approximations for nonexpansive mappings,
 {\it Bull. Amer. Math. Soc.}  \textbf{73} (1967), pp. 591--597.

\bibitem{Po} B.T. Polyak, Introduction to Optimization, Optimization Software, New York, 1987.

\bibitem{Rock}{ R.T. Rockafellar},  Monotone operators associated with saddle-functions and mini-max problems,
in {\em Nonlinear Functional Analysis,  Proceedings of Symposia in Pure Mathematics}, edited by F. Browder,
 American Mathematical Society  \textbf{18} (1), (1976), pp.~241--250.

\bibitem{Sch} M. Schatzman, A class of nonlinear differential equations of second order in time,
 {\it Nonlinear Analysis} \textbf{2} (1978), pp.~355--373.



\end{thebibliography}
\end{document}